\documentclass[11pt]{article}

\usepackage[dvipsnames]{xcolor}
\definecolor{SNSRed}{RGB}{183,42,58}
\definecolor{SNSBlue}{RGB}{0,114,142}
\definecolor{lilla}{RGB}{225,143,247}

\usepackage[colorlinks=true,hyperindex,linkcolor=SNSRed, pagebackref=false, citecolor=SNSBlue, urlcolor=SNSBlue, pdfpagelabels]{hyperref}
\usepackage{amsxtra}
\usepackage{amscd}
\usepackage{amsthm}
\usepackage{amsfonts}
\usepackage{amssymb}
\usepackage{eucal}
\usepackage{tikz-cd}
\usepackage[matrix,arrow,curve]{xy}
\usepackage{xy}
\usepackage[T1]{fontenc}
\usepackage[utf8]{inputenc}
\usepackage{mathtools}
\usepackage[margin = 3cm]{geometry}

\theoremstyle{plain}
\newtheorem{lem}{Lemma}[subsection]
\newtheorem{prop}[lem]{Proposition}
\newtheorem{cor}[lem]{Corollary}
\newtheorem{thm}[lem]{Theorem}

\theoremstyle{definition}
\newtheorem{deffa}[lem]{Definition}
\newtheorem{ex}[lem]{Example}
\newtheorem{rem}[lem]{Remark}

\theoremstyle{plain}
\newtheorem{thmsezione}{Theorem}[section]
\newtheorem{propsezione}[thmsezione]{Proposition}

\theoremstyle{definition}
\newtheorem{exsezione}[thmsezione]{Example}
\newtheorem{remsezione}[thmsezione]{Remark}
\newtheorem{conssezione}[thmsezione]{Construction}

\numberwithin{equation}{section}

\newcommand\Z{\mathbb{Z}}
\newcommand\Q{\mathbb{Q}}

\newcommand\calC{\mathcal{C}}
\newcommand\calF{\mathcal{F}}
\newcommand\calO{\mathcal{O}}
\newcommand\calI{\mathcal{I}}
\newcommand\calG{\mathcal{G}}
\newcommand\calX{\mathcal{X}}

\newcommand\frp{\mathfrak{p}}

\newcommand\pmeninf{p^{-\infty}}
\DeclareMathOperator{\ten}{\otimes}
\DeclareMathOperator{\im}{Im}

\newcommand{\Gm}{\mathbb{G}_m}
\newcommand{\Ga}{\mathbb{G}_a}
\newcommand{\Fp}{\mathbb{F}_p}
\newcommand{\Zp}{\mathbb{Z}_p}
\newcommand{\Qp}{\mathbb{Q}_p}
\newcommand{\mup}{\mu_{p}}                                              
\newcommand{\mupn}{\mu_{p^n}}

\DeclareMathOperator{\undis}{\coprod}
\newcommand{\onto}[2][]{\xrightarrow[#1]{#2}\mathrel{\mkern-14mu}\rightarrow}          
\newcommand\xto[2]{\xrightarrow[#1]{#2}}

\DeclareMathOperator{\Hom}{Hom}                                      
\DeclareMathOperator{\rk}{rank}                                        
\DeclareMathOperator{\coker}{coker}
\DeclareMathOperator{\op}{op}                                            
\DeclareMathOperator{\alg}{-Alg}                                        
\DeclareMathOperator{\sch}{Sch}                                         
                                         
\DeclareMathOperator{\D}{D}                                             
\DeclareMathOperator{\Tot}{Tot}                                         
\DeclareMathOperator{\Cone}{Cone}                                       
\DeclareMathOperator{\Ab}{\textit{Ab}}  
\DeclareMathOperator{\Alb}{\mathrm{Alb}}                                
\DeclareMathOperator{\Top}{\textit{Sh}}                                
\DeclareMathOperator{\Rgam}{\text{R}\Gamma}                            
\DeclareMathOperator{\Spec}{Spec}
\newcommand{\W}{W}               
\newcommand{\K}{K}
                                   
\newcommand{\Ll}{L}                                                     
\DeclareMathOperator{\Rlim}{Rlim}      
\DeclareMathOperator{\Rinf}{R\Gamma_\textnormal{inf}}                   
\DeclareMathOperator{\Rder}{\text{R}}
\DeclareMathOperator{\Rflat}{R\Gamma_{\textnormal{fppf}}}               
\DeclareMathOperator{\Ret}{R\Gamma_{\textnormal{\'{e}t}}}               
\DeclareMathOperator{\RdR}{R\Gamma_{\textnormal{dR}}}
                                         \DeclareMathOperator{\Br}{Br}                                           \DeclareMathOperator{\NS}{NS}
\newcommand\Hflat[1]{H^{#1}_{\textnormal{fppf}}}                        \newcommand\Het[1]{H^{#1}_{\textnormal{\'{e}t}}}                       
\DeclareMathOperator{\qcoh}{QCoh}
\newcommand\HdR[1]{H^{#1}_{\textnormal{dR}}}                            
\DeclareMathOperator{\Acr}{\mathbb{A}_{cris}}                           
\DeclareMathOperator{\Ainf}{\mathbb{A}_{inf}}                           
\newcommand\perf[1]{#1^{\flat}}                                         
\newcommand\coperf[1]{#1_{\textnormal{perf}}}                           
   
\DeclareMathOperator{\inff}{Inf}        
\DeclareMathOperator{\Cris}{CRIS}                                       
\DeclareMathOperator{\cris}{Cris}                                       
                                            
\DeclareMathOperator{\Rcris}{R\Gamma_{cris}}                            
\newcommand\ocris{\calO_{\textnormal{cris}}}
\newcommand\oinf{\calO_{\textnormal{inf}}}                           
\newcommand\icris{\calI_{\textnormal{cris}}}                            
\DeclareMathOperator{\Nyg}{F_N^1}                                      
\newcommand\Hcris[1]{H^{#1}_{\textnormal{cris}}}        
\newcommand\Hinf[1]{H^{#1}_{\textnormal{inf}}}                       
\DeclareMathOperator{\Spf}{Spf}

\renewcommand{\qed}{\hfill\square\vspace{.5em}}

\newcommand{\eu}{\begin{center}
\begin{minipage}{5cm}
\includegraphics[width=4.8cm]{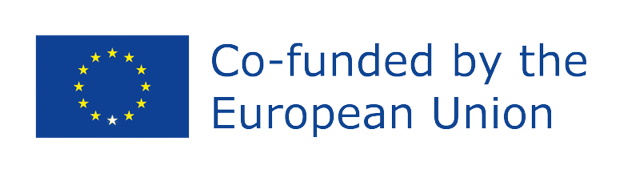}
\end{minipage}
\begin{minipage}{9.5cm}
\footnotesize Co-Funded by the European Union. Views and opinions expressed are however those of the author(s) only and do not necessarily reflect those of the European Union. Neither the European Union nor the granting authority can be held responsible for them.
\end{minipage}
\end{center}}

\begin{document}

\title{Some applications of the Nygaard filtration and quasisyntomic descent in positive characteristic}
\author{Livia Grammatica\footnote{Email: \textit{livia.grammatica@math.unistra.fr}\\
Institut de Recherche Math\'{e}matique Avanc\'{e}e (IRMA), Universit\'{e} de Strasbourg, 7 rue Ren\'{e} Descartes, 67000 Strasbourg, France.}} 
\date{\today}

\maketitle

\begin{abstract} This article gives an expository account of quasisyntomic descent and the first piece of the Nygaard filtration in positive characteristic, complemented by several applications to $p$-adic cohomology theories. The guiding result is a new approach to Illusie's comparison between fppf cohomology with $\Zp(1)$ coefficients and the slope $1$ part of crystalline cohomology. We follow work of Bhatt-Lurie, but give a more elementary presentation which does not rely on the formalism of $\infty$-categories. We then revisit Ogus' comparison theorem between infinitesimal cohomology and étale cohomology, and give new proofs of several results on fppf cohomology that were previously obtained with the de Rham-Witt complex. We also determine the action of multiplication-by-$n$ on the fppf cohomology of an abelian variety, answering a question of A. Skorobogatov to the author. This is an expanded version of the author's master thesis.
\end{abstract}

\tableofcontents
\section{Introduction}

Let $p$ be a prime number, $k$ a perfect field of characteristic $p$, $X$ a smooth and proper variety over $k$. Denote by $F$ the Frobenius acting on $\Hcris{i}(X/\W(k))$. We will be interested in the following comparison result, which is due to Illusie.

\begin{thmsezione}\label{mainthm}\textnormal{(\cite[Théorème II.5.5]{illusiedrw})}
If $k$ is algebraically closed, there is a natural map
\[
\Hflat{i}(X,\Zp(1))\to \Hcris{i}(X/\W(k))^{F=p}
\]
which becomes an isomorphism after inverting $p$\footnote{In this paper $\Hflat{i}(X,\Zp(1))$ denotes the $i$-th cohomology group of $\Rlim_n\Rflat(X,\mupn)$, but work of Milne shows that under our hypotheses the projective system of abelian groups $\{\Hflat{i}(X,\mupn)\}_n$ is Mittag-Leffler, so that $\Hflat{i}(X,\Zp(1))=\varprojlim_n\Hflat{i}(X,\mupn)$. See \cite[p.627]{illusiedrw} for further details.}.

\end{thmsezione}

The first goal of this paper is to present a proof of Theorem \ref{mainthm} using the modern approach to $p$-adic cohomology theories pioneered by Bhatt, Morrow, Scholze, Lurie, Drinfeld et al. To this effect we follow closely section $7$ of the paper of Bhatt-Lurie \cite{bhattlurie}. 

Illusie's original proof relies on the theory of the de Rham-Witt complex as the main tool for studying crystalline cohomology \cite{illusiedrw}. For a smooth $k$-scheme $X$, Illusie defines a projective system $(\W_n\Omega^{\bullet}_X)_{n\ge1}$ of complexes of sheaves over $X$ with the property that the hypercohomology of $W_n\Omega^{\bullet}_X$ computes the crystalline cohomology of $X$ over $\W_n(k)$. The proof of Theorem \ref{mainthm} is then essentially differential in nature, relying on properties of the Cartier isomorphism to relate fppf cohomology with logarithmic differential forms.

In Bhatt--Lurie's approach the focus is on the Nygaard filtration, and quasisyntomic descent is an essential tool to reduce statements to situations where cohomology can be described explicitly. Bhatt--Lurie's paper relies heavily on the formalism of $\infty$-categories, and part of the originality of this work is to present their proof in a more elementary way. We thus hope to make it accessible to geometers familiar with the basic tools of homological algebra (cohomology of sheaves on a site and derived categories) and crystalline cohomology.

In writing a detailed account of the techniques used in the proof we include a few applications of descent and the Nygaard filtration, some new and some well-known. As a direct application of quasisyntomic descent, we give a quick proof of a theorem of Ogus \cite{ogusinf} which identifies infinitesimal cohomology with the unit-root part of crystalline cohomology. Let $\Rinf(X/\W(k))$ denote the infinitesimal cohomology of $X$ (a definition is given in section \ref{infinitesima}).

\begin{thmsezione}\label{}\textnormal{(Theorem \ref{infmio})} If $k$ is a perfect field of characteristic $p$ and $X$ is a smooth $k$-scheme there is a natural isomorphism
\begin{equation}\label{isoinf}
\alpha_X:\Rinf(X/\W(k))\simeq\Rlim_F\Rcris(X/\W(k)),
\end{equation}
where $F$ is the Frobenius acting on $\Rcris(X/\W(k))$.

\end{thmsezione}

We then turn our attention to the long exact sequence
\[
\cdots\to{\Hflat{i}(X,\Zp(1))}\to{\Nyg\Hcris{i}(X/\W(k))}\xto{}{F/p-1}{\Hcris{i}(X/\W(k))}\to{\Hflat{i+1}(X,\Zp(1))}\to\cdots
\]
coming from triangle \eqref{triangolomainthm}. By studying the properties of the Nygaard filtration, we are able to prove many properties of fppf cohomology that Illusie obtains with the theory of the de Rham-Witt complex. 

We start with a general result proved by Illusie--Raynaud and Milne on the structure of fppf cohomology groups.

\begin{thmsezione}\label{}\textnormal{(Theorem \ref{finiteexponent})} For all $i\ge0$ the $\Zp$-module $\Hflat{i}(X,\Zp(1))$ is isomorphic to the direct sum of $\Hflat{i}(X,\Zp(1))_{\textnormal{tors}}$ with a free $\Zp$-module of rank $r=\rk_{\W(k)}\left(\Hcris{i}(X/\W(k))^{F=p}\right)$. Moreover, $\Hflat{i}(X,\Zp(1))_{\textnormal{tors}}$ is a $p$-group of finite $p$-exponent.
\end{thmsezione}

With different techniques one can get more precise results, namely that the group $\Hflat{i}(X,\Zp(1))_{\text{tors}}$ is the extension of a finite group by the $k$-points of a unipotent affine group scheme. This is proved in \cite{illusieraynaud} and \cite{milnezeta} (see Proposition\footnote{Note that in the statement one should replace ``finitely generated torsion $W(k)$-module''  with ``group of $k$-points of a commutative unipotent algebraic group''.} $2.2$ of \cite{geissermotivic} or Appendix A of \cite{skorkunneth}). 

We then turn our attention to the groups $\Hflat{i}(X,\Zp(1))$ for small values of $i$. For $i=1,2$ these groups turn out to be finite type $\Zp$-modules, but for $i=3$ the torsion subgroup can be infinite, even for abelian varieties. We give an example of this phenomenon at the end of the article, which hopefully serves to show that the formalism presented here is well-suited for doing explicit calculations. 
 
Finally, we identify a class of varieties for which fppf cohomology is completely determined by crystalline cohomology. These are essentially the varieties which satisfy the hypotheses of the Mazur "Newton-above-Hodge" theorem \cite[Theorem 8.26]{berthelotogus}. We call them straight varieties, but they are sometimes found in the literature under the name of Mazur--Ogus varieties \cite{ekedahl}. Some examples of straight varieties are abelian varieties, complete intersections, and K3 surfaces. We end our analysis by determining the action of the multiplication-by-$n$ map on the fppf cohomology of an abelian variety.

\begin{propsezione}\label{}\textnormal{(Corollary \ref{moltiplicazionen})} If $A$ is an abelian variety over an algebraically closed field $k$ of characteristic $p$, the multiplication-by-$n$ map $[n]:A\to A$ acts as $n^{i}$ on $\Hflat{i+1}(A,\Zp(1))_{\textnormal{tors}}$ and on $\Hflat{i}(A,\Zp(1))/\textnormal{tors}$.
\end{propsezione}

Part of the motivation for studying fppf cohomology is that the groups $\Hflat{i}(X,\Zp(1))$, for $i=2,3$, are closely related to the $p$-torsion of the Brauer group of $X$. Via the Kummer sequence, one can prove that there is a short exact sequence
\begin{equation}\label{brauer}
0\to\Br(X)[p^{\infty}]_{\mathrm{div}}\to\Br(X)[p^{\infty}]\to\Hflat{3}(X,\Zp(1))_{\text{tors}}\to0
\end{equation}
of abelian groups, where $\Br(X)[p^{\infty}]$ is the union of the $\Br(X)[p^n]$, and $\Br(X)[p^{\infty}]_{\mathrm{div}}$ is its maximal divisible subgroup. Furthermore, $\Br(X)[p^{\infty}]_{\mathrm{div}}$ is isomorphic to $\left(\Qp/\Zp\right)^{a}$ where $a$ is equal to $\dim_{\Qp}\Hflat{2}(X,\Qp(1))-\rk\NS(X)$, see \cite[Proof of Theorem A.1]{skorkunneth}. The upshot is then that $\Br(X)[p^{\infty}]$ is the direct sum of a divisible part and of a $p$-group of finite $p$-exponent, which comes from fppf cohomology.

From the above discussion, one sees that if $X$ is straight it is possible to write down a formula for the $p$-torsion of the Brauer group exclusively in terms of crystalline cohomology - so, for abelian varieties, exclusively in terms of its Dieudonné module. This will be exploited in a future paper with Alexei Skorobogatov and Yuan Yang.

\subsection*{Sketch of the proof of Theorem \ref{mainthm}} Let us briefly sketch the strategy of proof. If $X$ is a $k$-scheme, the Nygaard filtration is a filtration of crystalline cohomology $\Rcris(X/\W(k))$ in the derived category. For our purposes we only need the first piece of the Nygaard filtration, denoted $\Nyg{}\Rcris(X/\W(k))$, which by design sits in an exact triangle 
\[
\Nyg{}\Rcris(X/\W(k))\to\Rcris(X/\W(k))\to\Rgam(X,\calO_X)
\]
in $\D(\Zp)$, where the right-hand map is a natural augmentation map. Following Bhatt--Lurie we can construct a ``completed'' first Chern class map
\[
\hat{c}_1:\Rflat(X,\Zp(1))\to\Nyg{}\Rcris(X/\W(k))
\]
which is a refinement of the usual crystalline first Chern class. The composition
\[
\Rflat(X,\Zp(1))\xto{}{\hat{c}_1}\Nyg{}\Rcris(X/\W(k))\to\Rcris(X/\W(k))
\]
will be the comparison map of Theorem \ref{mainthm}. At the end of Section \ref{pruffa}, we will deduce Theorem \ref{mainthm} from the following result. Stated in this form, it goes back to the paper of Fontaine-Messing \cite{fontainemessing} where the syntomic topology is first introduced.

\begin{thmsezione}\label{mainthmbl}\textnormal{(\cite[Theorem 7.3.5]{bhattlurie}, Theorem \ref{triangololiscio})} If $X$ is a smooth scheme over a perfect field $k$, then
\begin{equation}\label{triangolomainthm}
\Rflat(X,\Zp(1))\xto{}{\hat{c}_1}\Nyg{}\Rcris(X/\W(k))\xto{}{F/p-1}\Rcris(X/\W(k))   
\end{equation}
is an exact triangle in $\D(\Zp)$.
\end{thmsezione}

The map $F/p-1$ is defined in Section \ref{pruffa}. In order to prove Theorem \ref{mainthmbl} we will use the technique of cohomological descent, which we briefly explain. The first step is to reduce to the case where $X=\Spec(R)$ is affine. The second step is to show that if $G$ is one of the three functors in \eqref{triangolomainthm}, $G$ satisfies descent along the map 
\[
\Spec(\coperf{R})\to\Spec(R),
\]
where $\coperf{R}=\varinjlim_{F_R}R$, which means that the natural map
\[
\xymatrix{G(\Spec(R))\ar[r]&\Tot\Bigl(G(\Spec(\coperf{R}))\ar[r]<1.5pt>\ar[r]<-1.5pt> & G(\Spec(\coperf{R}\otimes_R\coperf{R}))\ar[r]\ar[r]<-3pt>\ar[r]<3pt> & \dots\Bigr)}
\]
is an isomorphism in $\D(\Zp)$. By functoriality of triangle \eqref{triangolomainthm} this implies that it suffices to prove Theorem \ref{mainthmbl} when $X=\Spec(A)$, where $A$ is one of $\coperf{R},\coperf{R}\otimes_R\coperf{R},$ etc. 

The advantage of this reduction is that these rings are examples of what we call elementary quasiregular semiperfect (eqrsp) $\Fp$-algebras, which are introduced in Subsection \ref{elqrsp}. The point is that crystalline cohomology of an eqrsp algebra $A$ can be described very explicitly: it is isomorphic to $\Acr(A)[0]$, where $\Acr(A)$ be is a sort of ring of power series. The Nygaard filtration and fppf cohomology are also explicit and concentrated in degree zero so in the end we reduce to showing that a three-term sequence of abelian groups is exact. This is proved in the text in Theorem \ref{triangolosemipft}.

\subsection*{Outline of the article}

In Section \ref{prel} we present notation and results on crystalline cohomology which we will use throughout. We work with a version of the crystalline site which is slightly different from the usual, so we start by explaining the cohomology of sheaves in this setting. Then we explore in detail the crystalline cohomology of schemes of the form $\Spec(R)$, where $R$ is a semiperfect $\Fp$-algebra. 

In Section \ref{disc} we introduce a very light formalism for descent and show that crystalline cohomology satisfies descent along maps of the form $R\to\coperf{R}$ when $R$ is a Frobenius-smooth $\Fp$-algebra. This is a class of rings which includes smooth rings over a perfect field and more. We also give a proof of Ogus' theorem on infinitesimal cohomology.

In Section \ref{pruffa} we tackle the proof of the Theorem \ref{mainthmbl} and deduce Theorem \ref{mainthm}. Then, in the final section, we study the groups $\Hflat{i}(X,\Zp(1))$. We also determine $\Hflat{i}(X,\Zp(1))$ in two simple cases, to illustrate how these techniques can be used for computations. 

In the first four sections the only original result is the proof of Ogus' theorem. In Subsection \ref{semipft} we take many ideas from \cite{drinfeldacris}, and Section \ref{pruffa} is a rewriting of a portion of \cite[Section 7]{bhattlurie}. We have put effort in trying to present the results and arguments in an explicit and elementary way. Most of the results of Section \ref{applicazioni} are well-known but as far as we are aware the proofs presented here do not appear in previous literature.

\subsection*{Acknowledgements} This paper grew out of the author's masters thesis at Sorbonne Université, which was done under the supervision of Marco D'Addezio. The topic and most of the results presented here have been suggested by him, and we thank him heartily for his encouragements, his interest in this work and for many enlightening discussions. Thanks are due to Alexei Skorobogatov for suggesting to think about Corollary \ref{moltiplicazionen} and for several useful comments, and to Mauro Porta for explaining Remark \ref{mauro}. We also benefited from helpful conversations with Anna Cadoret, Yuan Yang, Emiliano Ambrosi and Haohao Liu. This project has received funding from the European Union’s Horizon Europe research and innovation programme under the Marie Skłodowska-Curie grant agreement n° 101126554.

\subsection*{Conventions and notations}
\begin{itemize}
    \item We fix the prime number $p$ throughout, as well as a perfect field $k$ of characteristic $p$. Unless otherwise stated, all rings and schemes are $\Fp$-algebras and $\Fp$-schemes respectively. All schemes are assumed to be quasi-compact and separated, for simplicity. 
    \item All divided power rings live over $(\Zp,(p))$, where $(p)$ is endowed with its unique divided power structure. All divided power envelopes are taken relative to $(\Zp,(p))$ 
    \item The $p$-typical Witt vectors are denoted by $\W(-)$, the Witt vectors of length $n$ by $W_n(-)$, and the Teichmüller representative of $a\in A$ by $[a]\in\W(A)$. The Frobenius automorphism of $\W(k)$ is $\sigma$. We write $\K=\W(k)[1/p]$, and for a $\W(k)$-module $M$ write $M[1/p]=M\otimes_{\W(k)}\K$. 
    \item If $A$ is a ring we write $\D(A)$ for the bounded below derived category -- note that this is not standard notation. 
    \item If $R$ is an $\Fp$-algebra we write $F_R$ for the Frobenius morphism $x\mapsto x^p$, simplified to $F$ when confusion is unlikely. The inclusion of perfect $\Fp$-algebras into $\Fp$-algebras has both a right and a left adjoint, respectively the perfection $R \mapsto \perf{R}=\varprojlim_{F}R$ and the coperfection $R \mapsto \coperf{R}=\varinjlim_{F}R$. Elements of $\perf{R}$ are sequences $(\dotsc,c_1,c_0)$ of elements of $R$ such that $c_i=c_{i+1}^p$. 
    \item Given a site $\calC$ the associated topos is $\Top(\calC)$ and the category of abelian sheaves is $\Ab(\calC)$. We denote the final object of $\Top(\calC)$ by $*$.
    \item A thickening of schemes is a closed immersion $X\to Y$ which induces a homeomorphism $|X|\to |Y|$. Equivalently, it is a closed immersion whose sheaf of ideals is locally nilpotent, meaning that every local section on an open affine subscheme is nilpotent.
    \item If $G$ is a $p$-divisible abelian group let $T_p(G)=\varprojlim_n G[p^n]$ denote the $p$-adic Tate module of $G$.
    \item An abelian group $A$ is a $p$-group if every element of $A$ is killed by a power of $p$. We say that $A$ has finite $p$-exponent if $A[p^n]=A$ for some $n$. The smallest such $n$, if it exists, is the $p$-exponent of $A$.
\end{itemize}

\section{Crystalline cohomology}\label{prel}


Crystalline cohomology was introduced by Grothendieck \cite{grothendieckinf} and developed by Berthelot \cite{berthelotcris} as a $p$-adic companion to $\ell$-adic cohomology. In this section we explain a simplified variant of Bhatt--Lurie's setup for working with crystalline cohomology. 

In Subsection \ref{convenzionicris} we define a big fppf version of the crystalline site. It is slightly different from Berthelot's definition but gives the same cohomology groups. We also gather some properties of crystalline cohomology that we will use in later proofs. The following sections cover some less classical topics in crystalline cohomology: derived $p$-completeness and the crystalline cohomology of semiperfect rings. In Subsection \ref{elqrsp} we introduce elementary quasiregular semiperfect algebras, which can be thought as contractible spaces for $p$-adic cohomologies. We work throughout in the derived category $\D(\Zp)$ where many constructions and statements become more natural.

\subsection{The big and small crystalline site}\label{convenzionicris}

The crystalline cohomology of an $\Fp$-scheme $X$ is typically defined as the cohomology of the structure sheaf on the small crystalline site. We briefly recall the definitions.

\begin{deffa}\label{}
(1) The \textit{small crystalline site} of $X$ is the category $\cris(X/\Zp)$ whose objects are divided power thickenings $(V,T,\delta)$, where $V$ an open subscheme of $X$. A morphism $(V',T',\delta')\to(V,T,\delta)$ is a morphism of divided power pairs where $V'\to V$ is a morphism of $X$-schemes. A family of maps $\{(V_i,T_i,\delta_i)\to(V,T,\delta)\}$ is a covering if $\{T_i\to T\}$ is a Zariski cover and $V_i=T_i\times_TV$ for all $i$. 

(2) The crystalline structure sheaf $\ocris$ is defined by $\ocris(V,T,\delta)=\Gamma(T,\calO_T)$. 

(3) The crystalline cohomology groups of $X$ are the cohomology groups of the complex $\Rgam(\cris(X/\Zp),\ocris)$, which we will abbreviate to $\Rcris(X/\Zp)$.

\end{deffa}

The small crystalline site is intimately linked to the Zariski topology on $X$. This means that other cohomology groups which are of interest to us, such as $\Hflat{n}(X,\Zp(1))$, can not be obtained as the cohomology of sheaves on this site. This limitation can be overcome by working with the more flexible big crystalline 
site, which we shall now define following \cite[Tag 07I5]{stacksproject}. While on the one hand it computes crystalline cohomology just like the small site, it has on the other hand many more sheaves to work with.

\begin{deffa}\label{}

(1) The \textit{big crystalline site} of $X$ is denoted by $\Cris(X/\Zp)$. Its objects are divided power
thickenings $(V,T,\delta)$, where $V$ is any $X$-scheme. A morphism $(V',T',\delta')\to(V,T,\delta)$ is a 
morphism of divided power pairs where $V'\to V$ is a morphism of $X$-schemes. A family of maps $\{(V_i,T_i,\delta_i)\to(V,T,\delta)\}$
is a covering if $\{T_i\to T\}$ is a covering for the fppf topology and $V_i=T_i\times_T V$.

(2) The structure sheaf is also denoted by $\ocris$, and is defined by $\ocris(V,T,\delta)=\Gamma(T,\calO_T)$. 

(3) Note that if we restrict the big crystalline site to an object $(V,T,\delta)$, we get the big fppf site of $T$. A sheaf $\calF$ is \textit{locally quasi-coherent} if it is quasi-coherent when restricted to any object $(V,T,\delta)$. For example, $\ocris$ is locally quasi-coherent, because $\ocris|_{(V,T,\delta)}=\calO_T$.

\end{deffa}

The following proposition echoes the fact that the fppf cohomology and the Zariski cohomology of a quasi-coherent sheaf coincide.

\begin{prop}\label{piattavzariski}

Let $\calF$ be a locally quasi-coherent sheaf of modules on $\Cris(X/\Zp)$. The natural map 
\[
\Rgam(\Cris(X/\Zp),\calF)\to\Rgam(\cris(X/\Zp),\calF|_{\cris(X/\Zp)})
\]
is a quasi-isomorphism. In particular, the cohomology of $\Rcris(X/\Zp)$ can be computed as the cohomology of $\Rgam(\Cris(X,\Zp),\ocris)$.

\end{prop}

\noindent{\it Proof.} By \cite[07IJ]{stacksproject} the inclusion $\cris(X/\Zp)\to\Cris(X/\Zp)$ induces a morphism of topoi $f:\Top(\Cris(X/\Zp))\to\Top(\cris(X/\Zp))$. Via the Leray spectral sequence of $f_{\*}$ we reduce to proving that $R^if_{*}\calF=0$ for $i>0$. For fixed $i$, the latter is the sheafification of $(V,T,\delta)\mapsto H^i((V,T,\delta),\calF|_{(V,T,\delta)})$. But this is zero whenever $V$ (or $T$) is affine (see \cite[Tag 07JJ]{stacksproject}). Since every object $(V,T,\delta)$ has a Zariski covering by affine divided power thickenings we find $R^if_{*}\calF=0$, which is what we wanted to prove.$\qed$

\begin{rem}\label{}

If $\calF$ is a sheaf on the big crystalline site we will sometimes write $\Rcris(X,\calF)$ for its cohomology.

\end{rem}

Let $\calF$ be a sheaf on the big fppf site of $X$. We can form a sheaf $\calF_c$ on $\Cris(X/\Zp)$ by setting $\calF_c(V,T,\delta)=\calF(V)$. If no confusion can arise we drop the subscript $c$ and write $\calF$ instead of $\calF_c$.

\begin{prop}\label{crisvfppf}

There is a canonical isomorphism $\Rcris(X,\calF_c)\simeq\Rflat(X,\calF)$. 

\end{prop}

\noindent{\it Proof.} Consider the morphism of topoi $u: \Top(\Cris(X/\Zp))\to\Top((X)_{\text{fppf}})$ induced by $(X)_{\text{fppf}}\to\Cris(X/\Zp)$ $V\mapsto(V,V,0)$. Since $u_*\calF_c=\calF$, by the Leray spectral sequence it is enough to show that $R^iu_{*}\calF_c=0$ for $i>0$. This sheaf is computed as the sheafification of $V\mapsto H^i((V,V,0),\calF_c|_{(V,V,0)})$. By definition, to give a covering of the object $(V,V,0)$ is equivalent to the give an fppf cover of $V$. Thus one has $H^i((V,V,0),\calF_c|_{(V,V,0)})=\Hflat{i}(V,\calF)$, and by locality of cohomology \cite[Tag 01FW]{stacksproject} we find that $R^iu_{*}\calF_c=0$. $\qed$

\begin{ex}\label{defnygaard}

Applying Proposition \ref{crisvfppf} to the fppf sheaf $\Ga$ we get an isomorphism $\Rcris(X,(\Ga)_c)\simeq\Rgam(X,\calO_X)$. There is a natural map of sheaves $\ocris\to(\Ga)_c$ which on $(V,T,\delta)$ corresponds to the surjection $\calO_T\onto{}{}\calO_V$, so in particular it is surjective. We call the kernel of this map $\icris$, and we shall denote by $\Nyg\Rcris(X/\Zp)$ the complex $\Rcris(X,\icris)$ and by $\Nyg H^n_{cris}(X/\Zp)$ its $n$-th cohomology group. It is called (the first piece of) the Nygaard filtration \cite[Remark 5.3.3]{bhattlurie}. By construction, there is an exact triangle
\begin{equation}
    \Nyg\Rcris(X/\Zp)\to\Rcris(X/\Zp)\to\Rgam(X,\calO_X)
\end{equation}
in the derived category. Since the complex $\Rgam(X,\calO_X)$ is $p$-torsion, we see that the maps $\Nyg\Rcris(X/\Zp)\to\Rcris(X/\Zp)$ and $\Nyg\Hcris{i}(X/\Zp)\to\Hcris{i}(X/\Zp)$ become isomorphisms after inverting $p$.

\end{ex}

\begin{rem}\label{notazionew}

When $X$ is a $k$-scheme one usually defines the sites $\cris(X/\W(k))$ and $\Cris(X/\W(k))$ as the divided power schemes $(V,T,\delta)$ such that $V$ an $X$-scheme and $T$ is a $\W(k)$-scheme. But in fact if $V$ is a $k$-scheme and $V\subseteq T$ is a divided power thickening then $T$ has a unique structure of $\W(k)$-scheme. This is proved in Lemma \ref{proprietauniversalewitt}. This justifies writing $\Rcris(X/\Zp)$ instead of $\Rcris(X/\W(k))$. When we want to emphasize the $\W(k)$-module structure of crystalline cohomology we will revert to writing $\Rcris(X/\W(k))$. 

\end{rem}

\begin{rem}
    It is also useful to define for all $n\ge0$ the truncated crystalline site $\Cris(X/\W_n(\Fp))$ whose objects are the $(V,T)\in\Cris(X/\Zp)$ such that $p^n=0$ in $T$, and then the truncated crystalline cohomology $\Rcris(X/\W_n(\Fp))$ is defined as the cohomology of the restriction of the structure sheaf to the truncated site. The isomorphism 
    \[
    \Gamma(\Cris(X/\Zp),\ocris)=\varprojlim_m\Gamma(\Cris(X/\W_m(\Fp)),\ocris)
    \]
    gives rise to a quasi-isomorphism of complexes
    \begin{equation}\label{limitecrisfinita}
    \Rcris(X/\Zp)=\Rlim_n\Rcris(X/\W_n(\Fp)).
    \end{equation}
    When $X$ is a $k$-scheme we may write $\Rcris(X/\W_n(k))$ instead of $\Rcris(X/\W_n(\Fp))$ to underline the $\W_n(k)$-module structure.
\end{rem}

\begin{rem}\label{}

When $X=\Spec(A)$ is an affine scheme we may write $\Cris(A/\Zp)$, $\Rcris(A/\Zp)$, etc. instead of $\Cris(X/\Zp)$, $\Rcris(X/\Zp)$, etc.

\end{rem}

We end this introductory section with a few important results on crystalline cohomology which we will need to use later on. Namely, we mention the comparison with de Rham cohomology, the perfectness of crystalline cohomology for proper schemes and the existence of a ``Verschiebung'' morphism.

\begin{prop}\label{confrontoderham}
Let $X$ be a smooth $k$-scheme. There is a canonical isomorphism
\[
\Rcris(X/\W(k))\otimes^L\Z/p\simeq\Rcris(X/\W_1(k))\simeq\RdR(X/k)
\]
\end{prop}

\noindent{\it Proof.}
See the results in \cite[Tag 07MI]{stacksproject}.
$\qed$

The first crystalline cohomology group $\Hcris{1}(X/\W(k))$ recovers a more classical object. 

\begin{prop}\label{crisabeliana}

Let $A$ be an abelian variety. There is a canonical isomorphism of $\Hcris{1}(A/\W(k))$ with the contravariant Dieudonné module of $A[p^{\infty}]$. Thus $\Hcris{1}(A/\W(k))$ is a free $\W(k)$-module of rank $2\dim A$.

\end{prop}

\noindent{\it Proof.}
This is a consequence of results of \cite{mazurmessing} as explained in the introduction of the book. Note that the authors work with a slightly different version of the crystalline site.
$\qed$

\begin{prop}\label{comparisonalbanese}
Let $X$ be smooth and proper, and let $a_X:X\to\Alb_X$ be a canonical morphism to the Albanese variety. The pull-back 
\[
a_X^*:\Hcris{1}(\Alb_X/\W(k))\to\Hcris{1}(X/\W(k))
\]
is an isomorphism. 
\end{prop}{}

\noindent{\it Proof.}
See \cite[II.3.11.2]{illusiedrw}. We haven't found in the literature a proof that doesn't use the de Rham-Witt complex.$\qed$

\begin{cor}\label{torsionehuno}

If $X$ is a smooth and proper $k$-scheme the $\W(k)$-module $\Hcris{1}(X/\W(k))$ is torsion-free.$\qed$

\end{cor}

We will only invoke Proposition \ref{comparisonalbanese} to use the Verschiebung $V$ on $\Hcris{1}(X)$, and to know that it acts compatibly with Verschiebung on de Rham cohomology, see Lemma \ref{diagrammacrisderhamlem}. The last fact we need is the finiteness of $\Hcris{i}(X/\W(k))$ when $X$ is proper.

\begin{prop}\label{cristallinaperfetta}
If $X$ is a smooth and proper $k$-scheme, the complex $\Rcris(X/\W(k))$ is a perfect object of $\D(\W(k))$. In particular, the cohomology groups $\Hcris{i}(X/\W(k))$ are $\W(k)$-modules of finite type.
\end{prop}

\noindent{\it Proof.}
See \cite[Tag 07MX]{stacksproject}.
$\qed$


\subsection{Derived $p$-completeness}\label{derivedpcomplete} Derived $p$-completeness is the appropriate analogue of $p$-completeness in the derived category $\D(\Z)$, and comes up naturally when working with $p$-adic cohomologies and derived categories. We only briefly discuss the definition and main results, and refer to existing literature for the proofs. The key points are the existence of the derived $p$-completion (Lemma \ref{derivedpcompletion}) and the derived Nakayama lemma (Lemma \ref{derivednakayama}).

\begin{deffa}\label{derivedpcomplete}

A complex $K\in\D(\Z)$ is $\textit{derived }p\textit{-complete}$ if the natural map $K\to\Rlim_n(K\ten^L_{\Z}\Z/p^n)$ is an isomorphism in $\D(\Z)$. The full subcategory of $\D(\Z)$ consisting of the derived $p$-complete objects will be denoted by $\D(\Zp)^{\wedge}$. It forms a triangulated subcategory of $\D(\Z)$.

\end{deffa}

This definition only partially agrees with the classical notion of $p$-completeness: if $M$ is an $\Zp$-module, then $M$ is $p$-complete in the usual sense if and only if $M[0]$ is derived $p$-complete and $\bigcap_np^nM=0$ (see e.g. \cite[Tag 091T]{stacksproject}).

\begin{lem}\label{equivderivedpcomplete}

Let $(K_n)$ be an inverse system of objects of $\D(\Z)$, and suppose that $p^n=0$ on $K_n$. The object $\Rlim_n\left(K_n\right)$ of $\D(\Z)$ is derived $p$-complete. In particular, for any $\Fp$-scheme $X$, $\Rcris(X/\Zp)$ is derived $p$-complete.

\end{lem}

\noindent{\it Proof.} The first assertion is \cite[Tag 091W]{stacksproject}. The second assertion follows from the first and from \eqref{limitecrisfinita}.
$\qed$

\begin{lem}\label{derivedpcompletion}

The inclusion $\D(\Zp)^{\wedge}\xhookrightarrow{}\D(\Z)$ has a left adjoint $K\mapsto K^{\wedge}$ given by the formula $K^{\wedge}=\Rlim_n(K\ten^L_{\Z}\Z/p^n)$. We call $K^{\wedge}$ the derived $p$-completion of $K$.

\end{lem}

\noindent{\it Proof.} This is \cite[Tag 0923]{stacksproject}.
$\qed$

\begin{ex}\label{completamentoGm}

If $X$ is a scheme we can consider the complex $\Rflat(X,\Gm)\in\D(\Z)$. Let us compute its derived $p$-completion:
\begin{align*}
    \Rflat(X,\Gm)^{\wedge}&=\Rlim_n\left(\Rflat(X,\Gm)\ten^L_{\Z}\Z/p^n\right)\\
    &=\Rlim_n\Cone\left(\Rflat(X,\Gm)\xto{}{\cdot p^n}\Rflat(X,\Gm)\right)\\
    &=\Rlim_n\Rflat(X,\mu_{p^n})[1]=\Rflat(X,\Zp(1))[1]
\end{align*}
Here we are writing $\Rflat(X,\Zp(1))$ for $\Rlim_n\Rflat(X,\mu_{p^n})$, as was explained in the footnote to the statement of Theorem \ref{mainthm}. 

Since tensoring with $\Z/p^n$ commutes with the derived limit, if $K\in\D(\Z)$ we have
\[
K^\wedge\otimes^L_{\Z}\Z/p^n\simeq\Rlim_m\big((K\ten^L_{\Z}\Z/p^m)\ten^L_{\Z}\Z/p^n\big)\simeq\Rlim_m\big((K\ten^L_{\Z}\Z/p^n)\ten^L_{\Z}\Z/p^m\big).
\]
The object $K\ten^L_{\Z}\Z/p^n$ is killed by $p^n$ and thus by Lemma \ref{equivderivedpcomplete} it is already derived $p$-complete. Hence $K^\wedge\otimes^L_{\Z}\Z/p^n\simeq K\otimes^L_{\Z}\Z/p^n$ and we get natural isomorphisms 
\[
\Rflat(X,\Zp(1))\otimes^L_{\Z}\Z/p^n\simeq\Rflat(X,\mupn{})
\]
for all $n\ge0$.
\end{ex}

We end this section with a result known as the derived Nakayama lemma.

\begin{lem}\label{derivednakayama}

Let $f:K_1\to K_2$ be a morphism in the derived $p$-complete category. If $\overline{f}:K_1\ten^L_{\Z}\Fp\to K_2\ten^L_{\Z}\Fp$ is an isomorphism then $f$ is an isomorphism.

\end{lem}

\noindent{\it Proof.} See \cite[Tag 0G1U]{stacksproject}.
$\qed$


\subsection{Semiperfect algebras and $\Acr$}\label{semipft}

An $\Fp$-algebra $B$ is \textit{semiperfect} if the Frobenius morphism $F_B$ is surjective, and is \textit{perfect} if $F_B$ is an isomorphism. The crystalline cohomology of semiperfect algebras enjoys a rather simple description. Before we come to it, let us recall the universal property of Witt vectors.

\begin{lem}\label{proprietauniversalewitt}

Let $B$ be a perfect $\Fp$-algebra and $S$ a $p$-complete ring. Any map $f:B\to S/(p)$ lifts uniquely to $\overline{f}:\W(B)\to S$.

\end{lem}

\noindent{\it Proof.} First, we lift $f$ to a multiplicative map $f_1:B\to S$. The strategy is to mimic the usual construction of Teichmüller representatives as in \cite[II, Proposition 8]{corpslocaux}: fix $b\in B$, and for every $n$ fix $y_n\in S$ such that $y_n$ reduces to $b^{1/p^n}$ mod $p$. The sequence $y_n^{p^n}$ converges to some element of $S$ independent of the choice of the $y_n$. If we set $f_1(b)=\lim_n\left(y_n^{p^n}\right)$ then $f_1$ is the desired multiplicative lift. Then, to define $\overline{f}$, write $x\in\W(B)$ as $x=\sum_n[b_n]p^n$ and set $\overline{f}(x)=\sum_nf_1(b_n)p^n$. One proves as in \cite[II, Proposition 10]{corpslocaux} that $\overline{f}$ is a ring homomorphism, and this settles the existence part. 

As for unicity, suppose another lift $g$ of $f$ is given. If $b\in B$ we have $g([b])=g([b^{1/p^n}])^{p^n}$, and since $g$ is a lift of $f$ the class of $g([b^{1/p^n]})$ mod $p$ is equal to $f([b^{1/p^n]})$. The argument above then shows that $g([b])=f_1([b])$, and finally by continuity we must have $g=\overline{f}$.
$\qed$

If $B$ is a semiperfect $\Fp$-algebra the natural map $B^\flat=\varprojlim_{F_B}B\to B$ is a surjection. The ring $\Acr(B)$ is defined as the divided power envelope of the surjection $\W(\perf{B})\onto{} \perf{B}\onto{q} B$.

\begin{prop}\label{Acris}

The divided power thickening $\Spec(B)\hookrightarrow\Spec(\Acr(B)/p^n\Acr(B))$ is the final object of $\Cris\left(B/\W_n(\Fp)\right)$.

\end{prop}

\noindent{\it Proof.} Let $(A,I)$ be an object in $\Cris\left(B/\W_n(\Fp)\right)$, so we have a ring map $f:B\to A/I$ and $p^nA=0$. By \cite[Tag 07GR]{stacksproject} there is an integer $N$ such that $x^{p^N}=0$ for all $x\in I$. We show that $f$ lifts uniquely to a map $\Acr(B)\to A$.

First, we claim that $f$ lifts uniquely to $f^{\flat}:\perf{B}\to A/p$: if $g,h$ are two such lifts then $g-h$ takes values in $I$, so if $x$ is an element of $\perf{B}$ we have $g(x)-h(x)=(g(x^{1/{p^N}})-h(x^{1/{p^N}}))^{p^N}=0$. Conversely, this gives a recipe for defining $f^{\flat}$ directly: if $x$ is in $\perf{B}$, let $a$ be a lift of $f(q(x^{1/{p^N}}))$ to $A/p$ and put $f^{\flat}(x)=a^{p^N}$. It is easy to check that this way we get a ring homomorphism lifting $f$.
    
Then, since $p^n=0$ in $A$, by Lemma \ref{proprietauniversalewitt} we get a unique lift of $f^{\flat}$ to a map $\W(\perf{B})\to A$. The universal property of the divided power envelope then induces the desired unique map $\Acr(B)\to A$ lifting $f$.$\qed$

\begin{rem}\label{ainf}

The ring $\W(\perf{B})$ is usually denoted by $\Ainf(B)$. In the proof of Proposition \ref{Acris} we have shown in passing that its spectrum represents the final object in the infinitesimal site of $B/\Zp$ (see Section \ref{infinitesima}).

\end{rem}

\begin{rem}
    In \cite[Appendix F]{bhattlurie}, the authors introduce and work with a ``$p$-completed'' variant of the big crystalline site. It has the advantage of containing objects such as $\Spec(B)\hookrightarrow\Spf(\Acr(B))$, which does not belong to our site, but it is also technically more demanding.
\end{rem}

\begin{cor}\label{cocrsp}

For a semiperfect algebra $B$ the complex $\Rcris(B/\Zp)$ is quasi-isomorphic to $\Acr(B)[0]$. More generally, if $\calF$ is a locally quasi-coherent sheaf of modules on $\cris(B/\Zp)$, and $M_n$ is its value on $\Acr/p^n\Acr$, then $\Rcris(B,\calF)$ is quasi-isomorphic to $\Rlim M_n$.$\qed$

\end{cor}


\subsection{Elementary quasiregular semiperfect algebras}\label{elqrsp}

Let $R$ be a smooth $k$-algebra. In Section \ref{disc} we show, loosely speaking, that the crystalline cohomology of $R$ can be recovered from the crystalline cohomology of the rings $\coperf{R}$,$\coperf{R}\ten_R\coperf{R}$, $\coperf{R}\ten_R\coperf{R}\ten_R\coperf{R}$ etc. These $\Fp$-algebras are all semiperfect, so their crystalline cohomology is governed by the behaviour of the functor $\Acr$.

The most important example for us is $R=A_n:=k[x_1,\dotsc,x_n]$. Setting $B_n=\coperf{(A_n)}$, we may write
\[
B_n=k[x_1^{\pmeninf},\dotsc,x_n^{\pmeninf}],
\]
and the tensor products $B_n\ten_{A_n}\dotsb\ten_{A_n}B_n$ are all isomorphic to rings of the form 
\[
B_n[t_1^{\pmeninf},\dotsc,t_m^{\pmeninf}]/(t_1,\dotsb,t_m),
\]
where $m=ns$, the integer $s$ being the number of factors in the product minus $1$.

\begin{deffa}\label{eqrsp}

An $\Fp$-algebra is \textit{elementary quasiregular semiperfect} (\textit{eqrsp} for short) if it is isomorphic to $B[x_1^{\pmeninf},\dotsc,x_m^{\pmeninf}]/(x_1,\dotsb,x_m)$ for some perfect $\Fp$-algebra $B$. An $\Fp$-scheme is \textit{locally eqrsp} if it admits a Zariski affine covering $\{\Spec(B_i)\}$ where the $B_i$ are eqrsp.

\end{deffa}

Eqrsp algebras are simple examples of quasiregular semiperfect algebras, see Definition \ref{spqrdef}. The point of this definition is that the $\Fp$-algebras $\coperf{R},\coperf{R}\ten_R\coperf{R},$ etc. are eqrsp up to Zariski localization.

\begin{prop}\label{modelloetale}

(1) For any prime ideal $\frp$ of $R$, there exists $f\in R\setminus\frp$ and étale morphisms 
\[
B_n\ten_{A_n}\dotsb\ten_{A_n}B_n\to(\coperf{R}\ten_R\dotsb\ten_R\coperf{R})_f.\]

(2) If $B$ is eqrsp and $B\to C$ is étale, then $C$ is eqrsp. In particular, all of the rings $(\coperf{R}\ten_R\dotsb\ten_R\coperf{R})_f$ are eqrsp.

\end{prop}

\noindent{\it Proof.} (1) Since $(\coperf{R})_f=\coperf{(R_f)}$ we may reduce to $R=R_f$ and suppose that there is an étale map $A_n\to R$. Then $(A_n\xto{}{F}A_n)\ten_{A_n} R$ is isomorphic to $R\xto{}{F} R$ (see \cite[Tag 0F6W]{stacksproject}) so we get a map
\[
B_n\to\coperf{R}=(A_n\to R)\ten_{A_n}B_n
\]
which is étale. The case of arbitrary tensor products follows from this one.

(2) We need to prove that if $A[x_1^{\pmeninf},\dotsc,x_n^{\pmeninf}]/(x_1,\dotsc,x_n)\to C$ is étale and $A$ is a perfect ring then $C$ is eqrsp. Let 
\[
q:A[x_1^{\pmeninf},\dotsc,x_n^{\pmeninf}]/(x_1,\dotsc,x_n)\onto{}A
\]
be the projection. If $I=\ker(q)$ the ring $S=C/IC$ is étale over $A$, and is therefore perfect. Now both $C$ and $S[x_1^{\pmeninf},\dotsc,x_n^{\pmeninf}]/(x_1,\dotsc,x_n)$ are étale over $A$ and reduce to $S$ mod $I$. Since $q$ is a thickening, by topological invariance of the étale site \cite[04DZ]{stacksproject} they are in fact isomorphic.$\qed$

We proceed to give a description of $\Acr$, and thus of crystalline cohomology, for eqrsp $\Fp$-algebras. This is done in \cite[section 2.5]{drinfeldacris} and we follow Drinfeld's exposition very closely. Write $C=B[x_1^{\pmeninf},\dotsc,x_n^{\pmeninf}]/(x_1,\dotsc,x_n)$ for some perfect algebra $B$ and set $I=(x_1,\dotsc,x_n)$ as an ideal of $B[x_1^{\pmeninf},\dotsc,x_n^{\pmeninf}]$. Set $\Z_+[1/p]=\{\alpha\in\Z[1/p]\text{ s.t. }\alpha\ge0\}$ and, if $\alpha=(\alpha_1,\dotsc,\alpha_n)\in\Z_+[1/p]^n$, write $x^\alpha=x_1^{\alpha_1}\dotsm x_n^{\alpha_n}$.

\begin{prop}\label{ainfspqr}

(1) Consider the map $B[x_1^{\pmeninf},\dotsc,x_n^{\pmeninf}]\to\perf{C}$ sending $x_i^{1/{p^n}}$ to $(\dotsc,x_i^{1/{p^{n+2}}},x_i^{1/{p^{n+1}}},x_i^{1/{p^n}})$. It induces an isomorphism between the $I$-adic completion of $B[x_1^{\pmeninf},\dotsc,x_n^{\pmeninf}]$ and $\perf{C}$, i.e. we can identify $\perf{C}$ with the ring of power series $\sum_{\alpha\in\Z_+[1/p]^n}a_{\alpha}x^{\alpha}$, $a_\alpha\in B$, such that for every $M>0$ the set $\{\alpha\le M\text{ s.t. }a_{\alpha}\ne 0\}$ is finite.

(2) The ring $\W(\perf{C})$ is identified with the ring of power series $\sum_{\alpha\in\Z_+[1/p]^n}b_{\alpha}x^{\alpha}$, $b_{\alpha}\in \W(B)$, such that for every $M>0$ the set $\{\alpha\le M\text{ s.t. }p^n\text{ does not divide }b_{\alpha}\}$ is finite.

\end{prop}

\noindent{\it Proof.} (1) We write the proof only for $n=1$, as the general case presents no significant difference. Set $x_1=x$. All elements in $I=\ker(\perf{C}\onto{}{}C)$ are of the form $(\dots,a_2,a_1,0)$, with $a_i\in C$ satisfying $a_i^{p^i}=0$. It is easy to see that then $a_i$ is divisible by $x^{1/p^{i}}$, so that $I$ is the principal ideal generated by $a=(\dotsc,x^{1/{p^2}},x^{1/{p}},0)$. The ring $\perf{C}$ is clearly $a$-complete, and we get a map $f:B^{\wedge}[x^{\pmeninf}]\to\perf{C}$ extending the map in the statement, where $B^{\wedge}[x^{\pmeninf}]$ denotes the $x$-completion of $B[x^{\pmeninf}]$. Both the domain and codomain of $f$ are $x$-complete and $x$-torsion-free, so $f$ is an isomorphism if $f$ mod $x$ is an isomorphism. But this is true, because $f$ reduces to the identity on $C$ modulo $x$.

(2) Let $R$ be the ring of formal series as in (2). By design $R$ is $p$-complete, $p$-torsion free and it reduces to $\perf{C}$ modulo $p$. Therefore $R\simeq\W(\perf{C})$.$\qed$

Note that $x^{1/p}$ can denote either the element $x^{1/p}\in C$ or its Teichmuller representative in $W(\perf{C})$. We will say explicitly which ring it belongs to whenever this can lead to confusion.

For a real number $y\ge0$ let $(y!)_p$ be the largest power of $p$ dividing $\lfloor y\rfloor!$, and for $\alpha=(\alpha_1,\dotsc,\alpha_n)\in\Z_+[1/p]^n$ set $(\alpha!)_p=(\alpha_1!)_p\dotsm(\alpha_n!)_p$. Let $F$ denote the Frobenius endomorphism of both $\W(\perf{C})$ and $\Acr(C)$. Following up on the notation of Example \ref{defnygaard} we write $\Nyg\Acr(C)$ for the kernel of $\Acr(C)\onto{}{}C$.

\begin{prop}\label{acrisspqr}

(1) We can identify $\Acr(C)$ with the ring of power series of the form 
\begin{equation*}
    \sum_{\alpha\in\Z_+[1/p]^n}b_{\alpha}\frac{x^{\alpha}}{(\alpha!)_p},\quad b_{\alpha}\in \W(B),
\end{equation*} 
such that for every $n>0$ the set $\{\alpha\text{ s.t. }p^n\text{ does not divide }b_{\alpha}\}$ is finite.

(2) The Frobenius $F$ acts as 
\[
F\left(\sum_{\alpha\in\Z_+[1/p]^n}b_{\alpha}\frac{x^{\alpha}}{(\alpha!)_p}\right)=\sum_{\alpha\in\Z_+[1/p]^n}F(b_{\alpha})\frac{((p\alpha)!)_p}{(\alpha!)_p}\frac{x^{p\alpha}}{((p\alpha)!)_p}.
\] 
Here $p\alpha=(p\alpha_1,\dots,p\alpha_n)$ and $\frac{((p\alpha)!)_p}{(\alpha!)_p}=p^{\lfloor\alpha\rfloor}=p^{\lfloor\alpha_1\rfloor+\dots+\lfloor\alpha_n\rfloor}$ is a nonnegative power of $p$.

(3) The ideal $\Nyg\Acr(C)$ is the set of power series $\sum_{\alpha\in\Z_+[1/p]^n}b_{\alpha}\frac{x^{\alpha}}{(\alpha!)_p}\in\Acr(C)$ such that $p$ divides $b_{\alpha}$ whenever $\alpha_i<1$ for all $i$.

\end{prop}

\noindent{\it Proof.} (1) The kernel of $\W(\perf{C})\onto{}{}C$ is generated by $p$ and the power series $x_i$ for $1\le i\le n$. Let $S$ be the smallest subalgebra of $\W(\perf{C})[1/p]$ containing $\W(\perf{C})$ and the $x^m/m!$ for $n\ge1$. It can be checked readily that $S$ satisfies the universal property of the divided power envelope of $\W(\perf{C})\onto{}{}C$, using that the $x_i$ are algebraically independent in $\W(\perf{C})$. Its $p$-completion is the ring of power series of (1).

(2) The Teichmuller representative $x^{1/p}\in W(\perf{C})$ satisfies $F(x^{1/p})=x$. Thus $F(x^{\alpha})=x^{p\alpha}$ and the rest follows easily.

(3) Point (1) tells us that $\Nyg\Acr(C)$ is the $p$-completion of the ideal of $S$ generated by $p$ and the $x_i^s/s!$ for $i=1,\dotsc,n$ and $s\ge1$, so (3) is another straightforward verification.$\qed$

\begin{cor}\label{filn}

The ideal $\Nyg\Acr(C)$ consists of the $a\in\Acr(C)$ for which $F(a)\in p\Acr(C)$.

\end{cor}

\noindent{\it Proof.} A power series $\sum_{\alpha\in\Z_+[1/p]^n}b_{\alpha}\frac{x^{\alpha}}{(\alpha!)_p}\in\Acr(C)$ is divisible by $p$ if and only if $p$ divides each one of its monomials, so we reduce to proving: $a=b_{\alpha}\frac{x^{\alpha}}{(\alpha!)_p}$ is in $\Nyg\Acr(C)$ if and only if $p$ divides $F(a)$. We have computed $F(a)=F(b_{\alpha})p^{\lfloor\alpha\rfloor}\frac{x^{p\alpha}}{((p\alpha)!)_p}$, and we know that $p$ divides $p^{\lfloor\alpha\rfloor}$ if and only if $\alpha_i\ge1$ for some $i$. The result follows from the description of $\Nyg\Acr(C)$ given in Proposition \ref{acrisspqr}.$\qed$

\begin{rem}\label{}

Let $D$ be a semiperfect algebra and call $S$ the divided power envelope of the surjection $W(\perf{D})\onto{}{}D$. The kernel of $S\onto{}{}D$ is generated by $p$ and the $[x]^n/n!$, where $n\ge1$ and $[x]$ is the Teichmüller representative of some $x\in\ker(\perf{D}\onto{}{}D)$. The Frobenius acts by sending $\frac{[x]^n}{n!}$ to $\frac{[x]^{pn}}{n!}=a\frac{[x]^{pn}}{(pn)!}$ where $a=\frac{(pn)!}{n!}$ is divisible by $p$, so it is always true that $F(\ker(\Acr(D)\onto{}{}D))\subseteq p\Acr(D)$. 

\end{rem}

In \cite{bhattlurie} the authors work with the more general notion of quasiregular semiperfect $\Fp$-algebras, which we now recall. The elementary quasiregular semiperfect algebras we consider are fundamental examples of quasiregular semiperfect algebras.

\begin{deffa}\label{spqrdef}

An $\Fp$-algebra $B$ is \textit{quasiregular semiperfect} (\textit{qrsp} for short) if it is semiperfect and the cotangent complex $\Ll_{B/\Fp}$ is a flat $B$-algebra supported in degree $-1$.

\end{deffa}

\begin{prop}\label{spqrprop}

(1) Being quasiregular semiperfect is an étale-local property.

(2) The quotient of a perfect $\Fp$-algebra by a Koszul-regular sequence is quasiregular semiperfect. In particular, elementary quasiregular semiperfect algebras are quasiregular semiperfect.

(3) The tensor products $\coperf{R}\ten_R\dotsb\ten_R\coperf{R}$ are quasiregular semiperfect.
\end{prop}

\noindent{\it Proof.} (1) If $B\to C$ is an étale cover, then $ \Ll_{C/\Fp}$ is quasi-isomorphic to $ \Ll_{B/\Fp}\ten_B^LC$, and 
\[
(B\xto{}{F}B)\ten_BC=(C\xto{}{F}C).
\]
Being surjective is an étale-local property for ring maps, and similarly for a complex to be flat and supported in degree $-1$. Thus $C$ is qrsp if and only if $B$ is qrsp.

(2) Such an algebra is semiperfect being a quotient of a perfect algebra. According to \cite[Tag 08SK]{stacksproject}, if $I$ is an ideal of $A$ generated by a Koszul-regular sequence, then $ \Ll_{(A/I)/A}$ is flat and concentrated in degree $-1$. If moreover $A$ is perfect then $ \Ll_{A/\Fp}=0$, because the Frobenius acts as $0$ and is an isomorphism, so $ \Ll_{(A/I)/A}= \Ll_{(A/I)/\Fp}$, and (2) follows.

(3) Follows from (1), (2) and Proposition \ref{modelloetale}.$\qed$

It follows that if $B$ is an eqrsp algebra, it is qrsp and $L_{B/\Fp}$ is perfect. The following example shows that these conditions are not enough to characterize eqrsp algebras.

\begin{prop}\label{}

Suppose $p\ne2$ and let $T=\Fp[x^{\pmeninf}]/(x^2)$. Then $T$ is qrsp, $L_{T/\Fp}$ is perfect but $T$ is not eqrsp.

\end{prop}

\noindent{\it Proof.}
The algebra $T$ is the quotient of a perfect algebra by a regular element, so $L_{T/\Fp}$ is perfect and by Proposition \ref{spqrprop} $T$ is qrsp. Suppose by contradiction that $T\simeq A[x_1^{\pmeninf},\dots,x_s^{\pmeninf}]/(x_1,\dots,x_s)$ for some perfect algebra $A$. Comparing the coperfection and the rank of the cotangent complex, we see that $A=\Fp$ and $s=1$. Thus there exists an isomorphism $g:\Fp[y^{\pmeninf}]/(y)\simeq T$.

Now note that every nonzero element of $T$ can be written uniquely in the form $ux^{i}$ where $0\le i<2$ is a rational number of the form $l/p^a$, and $u$ is a unit of $T$. Set $a_m=g(y^{1/p^m})$ and write $a_1=u_1x^{n}$ as above. Since $a_1^p=0$ we must have $2/p\le n$. Moreover, each of the $a_m$ is of the form $x^{n/p^{m-1}}u_m$ with $u_m$ a unit of $T$, and the $u_m$ satisfy $u_m=(u_{m+1})^p$.

By injectivity of $g$ we must have $(a_m)^{p^m-1}\ne0$. This equates to $n<2p^{m-1}/(p^{m}-1)$. Passing to the limit we get $n\le2/p$, so in the end $n=2/p$.

Any element of $k[y^{\pmeninf}]/(y)$ can be written uniquely in the form $vy^{i}$ where $0\le i<1$ is a rational number of the form $l/p^a$, and $v$ is a unit. Thus every element in the image of $g$ is of the form $ux^{2i}$ where $u$ is a unit of $T$. This implies that $x$ is not in the image of $g$, contradicting the assumption that $g$ is surjective. 

$\qed$

By mimicking the proof of Proposition \ref{modelloetale}, one shows that if $R\to S$ is étale then $S\simeq l[x^{\pmeninf}]/(x^2)$ for some étale $\Fp$-algebra $l$. A similar proof shows that $S$ is not eqrsp, so $R$ is not even étale-locally eqrsp.


\section{Cohomological descent}\label{disc}

Cohomological descent is a powerful extension of Grothendieck's descent theory for sheaves to derived categories. It is used to reduce constructions and statements about the cohomology of general spaces (schemes, varieties), to the same statement for a "nicer" class of spaces, where cohomology is easier to describe. Its origin can be traced back to Deligne's work, who uses descent to construct mixed Hodge structures for singular varieties, starting from Hodge structures of nonsingular ones, via iterated resolution of singularities. More recently it has been adopted by Scholze, Bhatt, Morrow, Lurie et al. as a foundational tool in $p$-adic cohomology theories. 

In Subsection \ref{generalitadiscesa} we present a very light formalism for cohomological descent, trading simplicity for flexibility. Then in Subsection \ref{frobeniussmooth} we study a special case of quasisyntomic descent for crystalline cohomology, namely descent along the map $\coperf{X}=\varprojlim_{F_X}X\to X$. The main result is Theorem \ref{discesafrobsmooth}. This will allow us to reduce statements about crystalline cohomology of smooth schemes to statements about $\Acr(R)$, where $R$ is an eqrsp algebra. The obvious advantage is that these rings are completely explicit. As a first application of this principle we give in Subsection \ref{infinitesima} a short proof of a result of Ogus comparing infinitesimal cohomology with the unit-root part of crystalline cohomology.

\subsection{Definition and examples}\label{generalitadiscesa}

Cohomological descent is typically an incarnation of the \v{C}ech spectral sequence for cohomology. We will formulate it in a way that makes maps between spectral sequences easier to visualise.

Let $\calC$ be a site and $K$ a sheaf of sets on $\calC$. The functor $\Ab(\calC)\to\Ab$ mapping $\calF$ to $\calF(K)=\text{Hom}_{\Top(\calC)}(K,\calF)$ is left exact, and we denote its right derived functor by $\Rgam(K,-)$. For example, if $K$ is the final object in $\Top(\calC)$ then $\Rgam(K,-)=\Rgam(\calC,-)$ is the usual derived global sections functor. Similarly, if $K$ is representable by some $X\in\calC$ then $\Rgam(K,-)=\Rgam(X,-)$.

\begin{prop}\label{discesacech}

Let $K'\to K$ be a surjective map of sheaves on $\calC$ and fix $\calF\in\Ab(\calC)$. Write $K'_m=K'\times_K\dotsb\times_KK'$ for the $(m+1)$-fold product of $K'$ over K. The natural map
\begin{equation}\label{totalization}
\xymatrix{\Rgam(K,\calF)\ar[r]&\Tot\Bigl(\Rgam(K',\calF)\ar[r]<1.5pt>\ar[r]<-1.5pt> & \Rgam(K'_1,\calF)\ar[r]\ar[r]<-3pt>\ar[r]<3pt> & \Rgam(K'_2,\calF)\dotsb\Bigr)}
\end{equation}
is an isomorphism in $\D(\Z)$.

\end{prop}

\begin{rem}\label{mauro}

Defining the right-hand term of \eqref{totalization} requires some care. The totalization of a cosimplicial object in \textit{cochain complexes} 
$$
\xymatrix{A_0\ar[r]<1.5pt>\ar[r]<-1.5pt> & A_1\ar[r]\ar[r]<-3pt>\ar[r]<3pt> & A_2\cdots}
$$
is defined as the totalization of the double complex that corresponds to it via the Dold-Kan correspondence. However taking the totalization of a cosimplicial object in the derived category is a much more subtle problem, as it is not clear that it can be lifted to an actual cosimplicial object in cochain complexes. The usual workarounds are working with the derived $\infty$-category, where a version of Dold-Kan holds (see Theorem 1.2.4.1 in \cite{higheralgebra}), or with simplicial topoi as in \cite[Tag 09VI]{stacksproject}.

Fortunately this will not be a problem in the situations we encounter: either $F(X)$ will be representable by a complex which is functorial in $X$, as is the case in Lemma \ref{discesaaciclica}, or $F=\Rgam(-,\calF)$ for some sheaf $\calF$ on a site. In the latter situation let $\calF\to\calI^{\bullet}$ be an injective resolution. Then
$$
\xymatrix{\Tot\Bigl(F(U)\ar[r]<1.5pt>\ar[r]<-1.5pt> & F(U\times_X U)\ar[r]\ar[r]<-3pt>\ar[r]<3pt> & F(U\times_X U\times_X U)\dotsb\Bigr)}
$$
can be identified with the totalization of the double complex $\calI^p(U^q)$ where $U^q=U\times_X\dotsb\times_XU$, $q+1$ times.

\end{rem}

\noindent{\it Proof of Proposition 3.1.1.} We need to check that the map induces an isomorphism upon taking cohomology. The cohomology of the left-hand side is simply $H^n(K,\calF)$. The complex on the right is the totalization of a double complex, so there is an $E_1$ spectral sequence converging to its cohomology. It is easy to see that $E_1^{p,q}=H^q(K'_p,\calF)$ and that the spectral sequence coincides with the spectral sequences of \cite[Tag 079Z]{stacksproject}, which converges to $H^{p+q}(K,\calF)$. Alternatively, we could apply \cite[Tag 0D8H]{stacksproject} with $E=\calF$ and the hypercovering being the \v{C}ech nerve of $K'\to K$.
$\qed$

Now we can give the definition of cohomological descent tailored to our purposes.

\begin{deffa}\label{}

Let $S$ be a scheme, $F:(\sch/S)^{\op}\to\D^{\ge0}(\Z)$ a functor, and let $f:U\to X$ be a map of $S$-schemes. We say that:

(1) $F$ satisfies descent along the map $f$ if the natural map
$$
\xymatrix{F(X)\ar[r]&\Tot\Bigl(F(U)\ar[r]<1.5pt>\ar[r]<-1.5pt> & F(U\times_X U)\ar[r]\ar[r]<-3pt>\ar[r]<3pt> & F(U\times_X U\times_X U)\dotsb\Bigr)}
$$
is an isomorphism in $\D(\Z)$. 

(2) $F$ satisfies descent for the étale, fppf, (quasi-)syntomic topology if it satisfies descent along any map
$f$ which is an étale, fppf, (quasi-)syntomic cover.

(3) $F$ satisfies descent for the Zariski topology if it satisfies descent for maps of the form $\undis U_i\to X$ where $\{U_i\}$ is a Zariski open covering of $X$. 

We give a similar definition in the affine case: let $A$ be a ring and $F:(A\alg)\to\D^{\ge0}(\Z)$ a functor.

(1) $F$ satisfies descent along the map of $A$-algebras $B\to C$ if the natural map
$$
\xymatrix{F(B)\ar[r]&\Tot\Bigl(F(C)\ar[r]<1.5pt>\ar[r]<-1.5pt> & F(C\ten_BC)\ar[r]\ar[r]<-3pt>\ar[r]<3pt> & F(C\ten_BC\ten_BC)\dotsb\Bigr)}
$$

is an isomorphism in $\D(\Z)$.

(2) $F$ satisfies descent for the étale, fppf, (quasi-)syntomic topology if it satisfies descent along any map of algebras which is an étale, fppf, (quasi-)syntomic etc. cover.

\end{deffa}{}

The typical blueprint for a descent argument is as follows: we wish to show $F$ satisfies some property on a certain class of schemes. We first show that $F$ satisfies descent for the Zarisi topology and reduce the problem to affine schemes in this class. Then, if $\Spec(R)$ is such a scheme, we look for a map $f:R\to R'$ such that for $F(R'),F(R'\ten_RR'),\dotsc$ the property becomes easy to prove. If $F$ satisfies descent along $f$, we can hope to deduce from this the property for $R$.

The simplest application of Proposition \ref{discesacech} is when $K,K'$ are representable sheaves and $K'\to K$ is a covering map in $\calC$.

\begin{cor}\label{discesagm}

Let $S$ be a scheme. The functor $\Ret(-,\Gm):(\sch/S)^{\op}\to\D(\Z)$ satisfies descent for the étale topology and for the fppf topology. 

\end{cor}{}

\noindent{\it Proof.} The first statement is a direct application of Proposition \ref{discesacech}. By a classical theorem of Grothendieck we have $\Ret(-,\Gm)=\Rflat(-,\Gm)$ (see \cite[Theorem 11.7]{brauer3}) so $\Ret(-,\Gm)$ also satisfies descent for the fppf topology.$\qed$

An interesting feature of this formalism is that it gives a way of saying that a sheaf is totally acyclic.

\begin{lem}\label{discesaaciclica}

Let $\calF$ be an abelian sheaf on a site $\calC$, and suppose that every covering of $\calC$ can be refined by a covering consisting of a single map (for example, $\calC$ is the small étale or fppf site over a scheme). If for every covering $Y\to X$ the map 
$$
\xymatrix{\calF(X)\ar[r]&\Tot\Bigl(\calF(Y)\ar[r]<1.5pt>\ar[r]<-1.5pt> & \calF(Y\times_X Y)\ar[r]\ar[r]<-3pt>\ar[r]<3pt> & \calF(Y\times_X Y\times_X Y)\dotsb\Bigr)}
$$
is an isomorphism in $\D(\Z)$, then for every object $X\in\calC$ we have $H^i(X,\calF)=0$ for $i>0$.

\end{lem}

\noindent{\it Proof.} Consider the the \v{C}ech spectral sequence 
\begin{equation*}
    E_1^{p,q}=H^q(Y_X^p,\calF)\implies H^{p+q}(X,\mathcal{F}),
\end{equation*}
where we write $Y_X^p$ for the $(p+1)$-fold product of $Y$ with itself over $X$. Since the bottom row of the spectral sequence is given by the \v{C}ech complex of $\mathcal{F}$ with respect to the cover $Y\to X$, our hypothesis says that $E_2^{p,0}$ is equal to $\calF(X)$ for $p=0$ and $0$ for $p>0$. By the short exact sequence 
\begin{equation*}
    0\to E^{1,0}_\infty\to H^1(X,\calF)\to E^{0,1}_\infty\subseteq H^1(Y,\calF)
\end{equation*}
we see that $H^1(X,\calF)$ injects into $H^1(Y,\calF)$, and this must hold for any covering $Y\to X$. But any cohomology class is locally trivial (see \cite[Tag 01FW]{stacksproject} for the precise statement), so $H^1(X,\calF)=0$. This holds for any $X\in\calC$.

If we consider the the filtration induced on $H^2(X,\calF)$ by the spectral sequence we now have that the first two terms are $0$, so that once again $H^2(X,\calF)$ injects into $H^2(Y,\calF)$ for any covering $Y\to X$. Arguing as above we see that $H^2(X,\calF)=0$ for all $X\in\calC$, and by induction one proves that all higher $H^i$s are $0$ as well.$\qed$

For crystalline cohomology we consider different sheaves for $K,K'$. Let $f:U\to X$ be a morphism of $\Fp$-schemes. The presheaf $h_U:(V,T,\gamma)\mapsto\Hom_X(V,U)$ on $\Cris(X/\Fp)$ is a sheaf for the fppf topology, and by definition we have $h_{U\times_X\dotsb\times_XU}=h_U\times\dotsb\times h_U$. Furthermore, by \cite[Tag 03F3]{stacksproject} if $\calF$ is a sheaf of abelian groups on $\Cris(X/\Zp)$ we have $\Rgam(h_U,\calF)=\Rgam(\Cris(U/\Zp),\calF|_U)$. Thus, if we want to show that the functor $\Rgam(\Cris(-/\Zp),\calF|_-)$ satisfies descent along $f$, it is enough to show that $h_U\to *$ is a surjective map of sheaves on $\Cris(X/\Zp)$, where $*$ denotes the final object of $\Cris(X/\Zp)$.

\begin{prop}\label{discesaetalecris}

The functor $\Rcris(-/\Zp):(\sch/\Fp)^{op}\to\D(\Zp)$ satisfies descent for the étale topology (and therefore also for the Zariski topology).

\end{prop}

\noindent{\it Proof.} We keep the notation from the previous paragraph, and we let $U\to X$ be an étale covering of $\Fp$-schemes. In the previous paragraph we saw that it suffices to show that $h_U\to *$ is surjective. In other words, we need to show that for an object $(V,T,\gamma)\in\Cris(X/\Zp)$ there is étale-locally a map to $h_U$. The projection $U\times_XV=V_U\to V\to V$ is an étale cover and $V\to T$ is a universal homeomorphism, so by topological invariance of the étale site \cite[Tag 04DZ]{stacksproject} there exists an étale $T$-scheme $T_U$ such that $T_U\times_T V=V_U$. Moreover, $T_U\to T$ is surjective, because $V_U\to V$ is. Call $\delta$ the unique PD-structure on $(V_U,T_U)$ compatible with that of $(V,T)$, which exists by \cite[Tag 07H1]{stacksproject}. Then $(V_U,T_U,\delta)\to(V,T,\gamma)$ is an étale cover with a map to $h_U$. This is what we wanted to prove. $\qed$


\subsection{Frobenius-smooth schemes}\label{frobeniussmooth}

Let $X$ be a smooth $k$-scheme. Our goal is to prove Theorem \ref{discesafrobsmooth}, which states that crystalline cohomology satisfies descent along the faithfully flat map $\coperf{X}\to X$. This holds for a larger class of schemes, called Frobenius-smooth in \cite{drinfeldstacky}, which we briefly review following Drinfeld.

\begin{deffa}\label{pbase}

(1) A \textit{finite} $p$-\textit{basis} of an $\Fp$-algebra $B$ is a finite set of elements $x_1,\dotsc,x_n$ such that each $b\in B$ can be written uniquely as a sum $\sum_{\alpha}b_{\alpha}^px^{\alpha}$, where $\alpha$ ranges over $\{0,\dotsc,p-1\}^n$ and $b_{\alpha}\in B$.

(2) A scheme $X$ is \textit{Frobenius-smooth} if every $x\in X$ has an affine neighbourhood $\Spec(B)$ such that $B$ has a finite $p$-basis. An $\Fp$-algebra $B$ is Frobenius-smooth if $\Spec(B)$ is.

\end{deffa}

Condition (1) says that $B$, considered as a $B$-module via the Frobenius, is free with basis $x_1^{\alpha_1}\dotsm x_n^{\alpha_n}$ where the $\alpha_i$ range over $\{0,\dotsc,p-1\}$. For completeness, we give a characterization of Frobenius-smooth schemes in terms of their Frobenius map. Recall that a morphism of schemes is \textit{syntomic} if it is flat, locally of finite presentation and the fibres are local complete intersections.

\begin{prop}\label{frobsintomico}

(\cite[Lemma 2.1.1]{drinfeldstacky}) The following are equivalent for an $\Fp$-scheme $X$:

(1) $X$ is Frobenius-smooth

(2) The Frobenius $F_X:X\to X$ is syntomic

\end{prop}

\noindent{\it Proof.} By the above observation, if (1) holds so does (2). Vice-versa, assume that (2) holds, and write $L_{X/X^p}$ for the cotangent complex of the morphism $F_X$. There is an exact triangle
$$
F_X^*\Ll_{X/\Fp}\to\Ll_{X/\Fp}\to\Ll_{X/X^p}
$$
in $\D(X)$ and the map $F_X^*\Ll_{X/\Fp}\to\Ll_{X/\Fp}$ is $0$ (\cite[Tag 0G5Z]{stacksproject}), so $L_{X/X^p}\simeq L_{X/\Fp}\oplus F_X^*L_{X/\Fp}[1]$. Because $F_X$ is syntomic, $L_{X/X^p}$ is a perfect complex with Tor-amplitude in $[-1,0]$ (\cite[Tag 08SL]{stacksproject}), which implies that both summands are perfect with Tor-amplitude in $[-1,0]$. Using that $F_X$ is faithfully flat, we see that $F_X^*L_{X/\Fp}[1]$ has the same Tor-amplitude as $L_{X/\Fp}[1]$ (\cite[Tag 068S]{stacksproject}), therefore both $L_{X/\Fp}$ and its shift $L_{X/\Fp}[1]$ have Tor-amplitude in $[-1,0]$. Thus $L_{X/\Fp}$ has Tor-amplitude $0$, so it is isomorphic to $H^0(L_X/\Fp)[0]\simeq\Omega_{X/\Fp}[0]$, and since it is perfect $\Omega_{X/\Fp}$ is locally free (\cite[Tag 0658]{stacksproject}). Then by \cite[Proposition 3.6]{pbasi} the scheme $X$ satisfies condition (1).$\qed$

Some examples of rings admitting a finite $p$-basis are the spectrum of a field of finite type over $\Fp$ and $\Spec(B[\![x_1,\dotsc,x_n]\!])$ for a perfect ring $B$. By \cite[Tag 0FW2]{stacksproject} the Frobenius morphism of a smooth scheme over a perfect ring is syntomic, so by Proposition \ref{frobsintomico} those schemes are also Frobenius-smooth.

In \cite{berthelotmessing} the authors study rings which have a non-necessarily finite $p$-basis. Using their results we now prove cohomological descent along the perfection map for Frobenius-smooth schemes.

\begin{thm}\label{discesafrobsmooth}

(1) Let $R$ be a Frobenius-smooth $\Fp$-algebra. The functors $\Rcris(-/\Zp)$ and $\Nyg{}\Rcris(-/\Zp)$ satisfy descent along the map $R\to\coperf{R}$.

(2) Let $X$ be a Frobenius-smooth $\Fp$-scheme. The functors $\Rcris(-/\Zp)$ and $\Nyg{}\Rcris(-/\Zp)$ satisfy descent along the map $\coperf{X}\to X$.

\end{thm}

\begin{rem}\label{}

Under this hypothesis, the Frobenius $F_R$ is faithfully flat and $\Spec(\coperf{R})\to\Spec(R)$ is a pro-fppf cover, but we do not expect crystalline cohomology to satisfy descent along all pro-fppf covers.

\end{rem}

\noindent{\it Proof.} (1) Recall that the functors $\Rcris(-/\Zp)$ and $\Nyg{}\Rcris(-/\Zp)$ are computed as the cohomology of the sheaves $\ocris$ and $\icris$, which are defined on the big crystalline site (see Example \ref{defnygaard} for the latter). Therefore, according to the discussion preceding Proposition \ref{discesaetalecris}, it is enough to show that $h_{\coperf{R}}\to*$ is a surjective map of sheaves on $\Cris(R/\Zp)$. This translates to showing that if $(\Spec(A/I)\hookrightarrow\Spec(A),\gamma)$ is an element of $\Cris(R/\Zp)$, there is a faithfully flat cover $(A,I,\gamma)\to(B,J,\delta)$ such that $R\to B/J$ factors through $R\to\coperf{R}$. We will take $J=IB$, and then by \cite[Tag 07H1]{stacksproject} $\delta$ exists and is uniquely determined. So we reduce to the following problem: find a faithfully flat $A$-algebra $B$ such that the composition $R\to A/I\to B/IB$ factors through $R\to\coperf{R}$.

We will use some deformation theory without proof. Namely, that for any $n$ there exists a lift $(R_n,F_n)$ of $(R,F_R)$ to $\W_n(k)$, or in other words there exists a flat $W_n(k)$-algebra $R_n$ and a flat ring map $F_n:R_n\to R_n$ such that $(R_n\xto{}{F_n}R_n)\ten k=R\xto{}{F_R}R$. An explicit construction of the lift is given in \cite[Corollary 1.2.7]{berthelotmessing}.

By definition of the crystalline site every element of $I$ is nilpotent, so we may choose $n$ such that $p^n=0$ in $A$. We claim that $R\to A/I$ lifts to $R_n\to A$. If $R$ is smooth this boils down to a version of Hensel's lemma, because $I$ is locally nilpotent. The general case follows from \cite[Proposition 1.2.6]{berthelotmessing}. Set $(\coperf{R})_n=\varinjlim_{F_n}(R_n)$. This is a lift of $\coperf{R}$ to $\W_n(\Fp)$ and $R_n\to(\coperf{R})_n$ is faithfully flat. Then $B=A\ten_{R_n}(\coperf{R})_n$ is the $A$-algebra we were looking for: is faithfully flat over $A$ and $B/IB=A/I\ten_{R}\coperf{R}$ so $R\to B/IB$ factors through $R\to\coperf{R}$.

(2) Let $\{U_i\to X\}$ be a finite Zariski cover, and set $U=\undis U_i$. Then $h_U\to*$ is surjective by Proposition \ref{discesaetalecris} and $h_{\coperf{U}}\to h_U$ is surjective by the proof of (1). Thus $h_{\coperf{U}}\to h_{\coperf{X}}\to*$ is surjective and $h_{\coperf{X}}\to*$ must be as well. By the discussion preceding the proof of Proposition \ref{discesaetalecris} this is enough to conclude.
$\qed$

\begin{rem}\label{}

The algebra in the proof of Theorem \ref{discesafrobsmooth} becomes clearer if $R=k[x_1,\dotsc,x_n]$. In that case, the map $R\to A/I$ corresponds to $n$ elements $a_1,\dotsc,a_n\in A/I$, and we seek a faithfully flat extension of $A$ whose reduction mod $I$ has all $p$-th power roots of the $a_i$. If we take lifts $b_i\in A$ of the $a_i$ we get a lift $\W_n(k)[x_1,\dots,x_n]\to A$ of $R\to A/I$, and then we can adjoin freely all $p$-th roots of the $b_i$ by letting $B=A\ten_{\W_n(k)[x_1,\dotsc,x_n]}\W_n(k)[x_1^{\pmeninf},\dots,x_n^{\pmeninf}]$. 

\end{rem}

\begin{cor}\label{discesaesplicita}

Let $R$ be a smooth $k$-algebra. The maps
\[
\xymatrix{\Rcris(R/\Zp)\ar[r]&\Tot\Bigl(\Acr(\coperf{R})\ar[r]<1.5pt>\ar[r]<-1.5pt> & \Acr(\coperf{R}\ten_R\coperf{R})\ar[r]\ar[r]<-3pt>\ar[r]<3pt> & \dotsb\Bigr)}
\]
and
\[
\xymatrix{\Nyg\Rcris(R/\Zp)\ar[r]&\Tot\Bigl(\Nyg\Acr(\coperf{R})\ar[r]<1.5pt>\ar[r]<-1.5pt> & \Nyg\Acr(\coperf{R}\ten_R\coperf{R})\ar[r]\ar[r]<-3pt>\ar[r]<3pt> & \dotsb\Bigr)}
\]
are quasi-isomorphisms. Similar formulas hold for smooth $k$-schemes.
$\qed$
\end{cor}

\begin{rem}\label{}

Let us mention a result of Drinfeld \cite[Theorem 2.4.2]{drinfeldstacky} which fits in naturally with the discussion above.

Denote by $\calC$ be the fibred category over $\Cris(X/\Zp)^{\text{op}}$ whose fibre over $(V,T,\delta)$ is $\qcoh(T)$, and morphisms $\calF'\to\calF$ over $f:(V',T',\delta')\to(V,T,\delta)$ are given by morphisms of quasi-coherent sheaves $\calF\to f_*(\calF')$. The category $\calC$ is a stack over the big crystalline site: this statement is no more and no less the fact that the fibred category of quasi-coherent sheaves over schemes form a stack for the fppf topology, see \cite[Theorem 4.23]{vistolidescent}.

This means that we can define $\calC(h_U)$ for any scheme $U$ over $X$. Indeed, if we choose a surjective map of sheaves $\undis U_i\onto{}h_U$, where the $U_i$ are objects of $\Cris(X/\Zp)$, we can define $\calC(h_U)$ as an equalizer:
$$
\xymatrix{\calC(h_U)\ar[r]&\prod_i\calC(U_i)\ar[r]<1.5pt>\ar[r]<-1.5pt> & \prod_{i,j}\calC(U_i\times_{h_U} U_j)\ar[r]\ar[r]<-3pt>\ar[r]<3pt> & \prod_{i,j,k}\calC(U_i\times_{h_U} U_j\times_{h_U} U_k)}
$$
The products appearing in the sequence are all objects of $\Cris(X/\Zp)$, so this definition makes sense. Note however that this is somewhat informal: to be rigorous one should understand this equalizer as a suitable category of descent data, as explained in \cite{vistolidescent}.  The point of this construction is that $\calC(h_U)$ is equivalent to the 
category of crystals over $U$.

Now, suppose that $X$ is Frobenius-smooth. We have proved that $h_{\coperf{X}}\to*$ is surjective in $\Cris(X/\Zp)$, so by the stack property we have a similar equalizer diagram
$$
\xymatrix{\calC(*)\ar[r]&\calC(h_{\coperf{X}})\ar[r]<1.5pt>\ar[r]<-1.5pt> & \calC(h_{\coperf{X}\times_X\coperf{X}})\ar[r]\ar[r]<-3pt>\ar[r]<3pt> & \calC(h_{\coperf{X}\times_X\coperf{X}\times_X\coperf{X}})}
$$

A crystal on a semiperfect scheme $Y$ is the same as a quasi-coherent sheaf on the $p$-adic formal scheme $\Acr(Y)$. Thus we get an equalizer diagram
$$
\xymatrixrowsep{-1pt}\xymatrixcolsep{3pt}\xymatrix{\calC(*)\ar[rr]& &\qcoh(W(\coperf{X}))\ar[rr]<1.5pt>\ar[rr]<-1.5pt> & & \qcoh(\Acr(\coperf{X}\times_X\coperf{X}))\\
&&\ar[r]\ar[r]<-3pt>\ar[r]<3pt> &&\qcoh(\Acr(\coperf{X}\times_X\coperf{X}\times_X\coperf{X}))}
$$

Let $\calG$ the divided power completion of the subscheme $\coperf{X}\times_X\coperf{X}\subseteq\W(\coperf{X})\times\W(\coperf{X})$. In the above cited theorem Drinfeld shows that
$$
\xymatrix{\dotsb\Acr(\coperf{X}\times_X\coperf{X}\times_X\coperf{X})\ar[r]\ar[r]<3pt>\ar[r]<-3pt>&\Acr(\coperf{X}\times_X\coperf{X})\ar[r]<1.5pt>\ar[r]<-1.5pt>&\Acr(\coperf{X})}
$$
is the \v{C}ech nerve of $\calG$ acting on $W(\coperf{X})$, and deduces that $\calG$ is a flat and affine groupoid acting on $\Acr(\coperf{X})=\W(\coperf{X})$. The consequence of this observation is that pulling back crystals induces an equivalence of categories
$$
\calC(*)\simeq\qcoh\left(\left[W(\coperf{X})/\calG\right]\right)
$$
Moreover, our proof of Theorem \ref{discesafrobsmooth} shows that if $\calF$ is a crystal on $X$ and $F$ is its pullback to $\calC(h_{\coperf{X}})=\qcoh(\W(\coperf{X}))$ then the natural map
$$
\Rgam(X,\calF)\to\Rgam([W(\coperf{X})/\calG],F)
$$
is an isomorphism. Indeed, coherent cohomology of a sheaf $M$ on a quotient stack $\calX=[A/H]$, with $H$ a flat and affine groupoid acting on $A$, is computed by the totalization of the cosimplicial complex
$$
\xymatrix{\Rgam(A,M|_{A})\ar[r]<-1.5pt>\ar[r]<1.5pt>&\Rgam(A\times_{\calX}A,M|_{A\times_{\calX}A})\ar[r]\ar[r]<3pt>\ar[r]<-3pt>&\Rgam(A\times_{\calX}A\times_{\calX}A,M|_{A\times_{\calX}A\times_{\calX}A})\dotsb}
$$
This is again just the \v{C}ech spectral sequence for the covering $A\to\calX$.

\end{rem}

We end this section with an enhancement of Corollary \ref{discesagm} which we will need in Section \ref{pruffa}.

\begin{lem}\label{discesaperfgm}

Let $X$ be Frobenius-smooth. The functor $\Ret(-,\Gm)$ satisfies descent along the map $\coperf{X}\to X$.

\end{lem}

\noindent{\it Proof.} If $Y=\varprojlim_n Y_n$ where $(Y_n)_{n\ge1}$ is a projective limit of quasi-compact quasi-separated schemes with affine transition maps, then the natural map $\varinjlim_n\Ret(Y_n,\Gm)\to\Ret(Y,\Gm)$ is an isomorphism in $\D(\Z)$, see e.g. \cite[Tag 09YQ]{stacksproject} for a proof. The conclusion follows from Corollary \ref{discesagm} and the fact that filtered colimits commute with totalization. $\qed$

\subsection{Infinitesimal cohomology in positive characteristic}\label{infinitesima}
Infinitesimal cohomology was introduced by Grothendieck \cite{grothendieckinf} as a first attempt at defining a $p$-adic Weil cohomology theory. We briefly recall the definition.

Let $S$ be any scheme, not necessarily over $\Fp$. If $X$ is an $S$-scheme, the objects of the (small) infinitesimal site $\inff(X/S)$ are pairs of $S$-schemes $(U,T)$ where $U$ is an open subscheme of $X$ and $U\subseteq T$ is a thickening with nilpotent ideal. A morphism $(U,T)\to (U',T')$ is an $S$-morphism of pairs where $U\to U'$ is an $X$-morphism (i.e. an open immersion). The infinitesimal site comes equipped with the structure sheaf $\oinf$ defined by $\oinf(U,T)=\calO_T(T)$. The infinitesimal cohomology of $X$ is the complex $\Rinf(X/S):=\Rgam\left(\inff(X/S),\oinf\right)$ and its cohomology groups are denoted by $\Hinf{i}(X/S)$.

A theorem of Grothendieck asserts that in characteristic zero infinitesimal cohomology coincides with de Rham cohomology.

\begin{thm}\label{infinitesimacarzero}(\cite[Theorem 4.1]{grothendieckinf}) If $S$ is a $\Q$-scheme, and $X$ is a smooth $S$-scheme, then $\Hinf{i}(X)$ is canonically isomorphic to $\HdR{i}(X/S)$.
\end{thm}

This result is remarkable as it gives a definition of de Rham cohomology without differential forms. In positive characteristic Theorem \ref{infinitesimacarzero} does not hold anymore, and infinitesimal cohomology lacks many properties of ``nice'' cohomology theories such as a cycle class map. Ogus described it when $S=\Spec(\W(k))$, where $k$ is a perfect field of characteristic $p$. Recall that if $M$ is a finite type $\W(k)$-module with a $\sigma$-linear endomorphism $F$, the unit-root part of $M$, denoted $M^0$, is the largest $F$-stable submodule on which $F$ is bijective. It can be identified with $\bigcap_nF^n(M)$.

\begin{thm}\label{infogus}(\cite[Theorem 4.9]{ogusinf})
Let $k$ be a perfect field of characteristic $p$ and let $X$ be a proper smooth $k$-scheme. There is a natural isomorphism
\[
\Hinf{i}(X/\W(k))\simeq\Hcris{i}(X/\W(k))^0
\]
for all $i$, which respects the action of Frobenius.
\end{thm}

The map comes from a natural construction: there is a functor from the infinitesimal site of $X$ to the small crystalline site of $X$, which maps $(U,T)$ to the divided power envelope of the thickening $U\subseteq T$, and induces morphisms $\Hinf{i}(X/\W(k))\to\Hcris{i}(X/\W(k))$ for all $i$. Furthermore, if $k$ is algebraically closed the unit-root part of $\Hcris{i}(X/\W(k))$ can be identified with $\Het{i}(X,\Zp)\otimes\W(k)$. This fact can be proved by descent (see Remark \ref{comparisonetale}) and Theorem \ref{infogus} thus establishes an isomorphism between infinitesimal cohomology and étale cohomology with $\W(k)$ coefficients. We propose the following restatement of Ogus' theorem.

\begin{thm}\label{infmio}

Let $k$ be a perfect field of characteristic $p$ and $X$ a smooth scheme. There is a natural isomorphism
\begin{equation}\label{isoinf}
\alpha_X:\Rinf(X/\W(k))\simeq\Rlim_F\Rcris(X/\W(k))
\end{equation}
in $\D(\Zp)$, where $F$ is the Frobenius acting on $\Rcris(X/\W(k))$.

\end{thm}

We prove Theorem \ref{infmio} via descent along the map $\coperf{X}\to X$. The first step is to consider the big fppf infinitesimal site as in Subsection \ref{convenzionicris}. The properties of the construction are identical to the crystalline case, so we omit them. Then we have the following descent properties for $\Rinf(-/\W(k))$.

\begin{prop}\label{discesainf}

Let $X$ be a smooth $k$-scheme.
\begin{enumerate}
    \item The functor $\Rinf(-/\W(k))$ satisfies descent for the étale topology.
    \item The functor $\Rinf(-/\W(k))$ satisfies descent along the map $\coperf{X}\to X$.
\end{enumerate}
\end{prop}

\noindent{\it Proof.} Point (1) is proved the same as Proposition \ref{discesaetalecris}, and point (2) the same as Theorem \ref{discesafrobsmooth}.
$\qed$

\begin{lem}\label{infainf}
    Let $R$ be a semiperfect $k$-algebra. The complex $\Rinf(\Spec(R)/\W(k))$ is quasi-isomorphic to $W(R^\flat)[0]=\Ainf(R)[0]$.
\end{lem}
\noindent{\it Proof.} As mentioned in Remark \ref{ainf}, the proof of  Corollary \ref{cocrsp} shows that $\Spec(\W_n(\perf{R}))$ is the initial object in the (truncated) big infinitesimal site. Then it remains to verify that cohomology can be computed in the big site, just like for crystalline cohomology
$\qed$

Now we can construct $\alpha_X$. Consider the commutative diagram
\[\begin{tikzcd}
	{\Rlim_F\Rinf(X/\W(k))} & {\Rlim_F\Rcris(X/\W(k))} \\
	{\Rinf(X/\W(k))} & {\Rcris(X/\W(k))}
	\arrow["a", from=1-1, to=1-2]
	\arrow["q", from=1-1, to=2-1]
	\arrow[from=1-2, to=2-2]
	\arrow[from=2-1, to=2-2]
\end{tikzcd}\]
We claim that $q$ is a quasi-isomorphism, and in fact that $F$ acts as a quasi-isomorphism on $\Rinf(X/\W(k))$. Arguing by descent (Proposition \ref{discesainf}) one can first reduce to a smooth affine $X$, and then to $X=\Spec(R)$ with $R$ eqrsp, where it follows from Lemma \ref{infainf}. Thus we can define $\alpha_X$ as the composition $a\circ q^{-1}$.

\vspace{.5em}
\noindent{\it Proof of Theorem 3.3.3.} Arguing once more by descent we reduce to the case where $X=\Spec(R)$ with $R$ eqrsp, where we are left to show that the map
\[
\Ainf(R)\to\bigcap F^n\left(\Acr(R)\right),
\]
which is an inclusion, is an isomorphism. This is straightforward to check using the explicit descriptions of $\Ainf(R)$ and $\Acr(R)$ given in Proposition \ref{ainfspqr} and Proposition \ref{acrisspqr}: let $a=\sum b_{\alpha}\frac{x^{\alpha}}{(\alpha!)_p}$ be an element lying in $F^n\Acr(R)$ for all $n$. Since the $\frac{x^{\alpha}}{(\alpha!)_p}$ form a topological basis of $\Acr(R)$ and $F$ preserves $\langle\frac{x^{\alpha}}{(\alpha!)_p}\rangle_{\W(B)}$, this implies that each $b_{\alpha}\frac{x^{\alpha}}{(\alpha!)_p}$ lies in this intersection. Writing 
\begin{equation}\label{equazioneinf}
    b_{\alpha}\frac{x^{\alpha}}{(\alpha!)_p}=F^n(c_n),
\end{equation}
it is clear that $c_n$ is of the form $b_n\frac{x^{\beta}}{(\beta!)_p}$. Then, applying inductively the formula for Frobenius from point (2) of Proposition \ref{acrisspqr}, equation \eqref{equazioneinf} implies that $\beta=\alpha/p^n$ and 
\begin{equation*}
    b_\alpha=\sigma^n(b_n)p^{\sum_{0\le i\le{n-1}}\lfloor p^i\beta\rfloor}=\sigma^n(b_n)p^{\sum_{1\le i\le{n}}\lfloor \alpha/p^{i}\rfloor}=\sigma^n(b_n)p^{\sum_{1\le i\le{n}}\lfloor\lfloor\alpha\rfloor/p^{i}\rfloor}.
\end{equation*}
By Legendre's formula, such a $b_n$ exists for all $n$ if and only if  $(\alpha!)_p$ divides $b_{\alpha}$, which is equivalent to $b_{\alpha}\frac{x^{\alpha}}{(\alpha!)_p}\in\Ainf(R)$.
$\qed$

\noindent{\it Proof of Theorem 3.3.2.} Now $k$ is perfect and $X$ is a smooth and proper $k$-scheme. Let $M$ be a finitely generated $\W(k)$-module with a Frobenius-linear endomorphism $F$. We claim that $\Rder^1\lim_FM$ is zero, and that $H^0(\Rlim_FM)$ is the unit-root part of $M$ concentrated in degree zero. We have 
\[
\Rlim_FM\simeq\Rlim_F\Rlim_n\left(M/p^nM\right)\simeq\Rlim_n\Rlim_F(M/p^nM).
\]
Since $M/p^nM$ is of finite $\W(k)$-length, the projective system $(M/p^nM,F)$ is Mittag-Leffler and the $\Rder^1\lim$ vanishes. The subgroups $F^i(M/p^nM)$ form a decreasing chain of $\W(k)$-submodules of $M/p^nM$ which must stabilise, so that $(M/p^nM)^0=F^N(M/p^nM)$ for $N>\!>0$. Therefore $\varprojlim_F(M/p^nM)\simeq(M/p^nM)^0$. 

We claim that $(M/p^nM)^0=(M^0+p^nM)/p^nM$. Indeed, the second module embeds naturally in the first module, and the quotient is a submodule of $M/(M^0+p^nM)$. By definition of $M^0$ the endomorphism $F$ is nilpotent on this space. Since $F$ acts bijectively on any quotient of $(M/p^nM)^0$ by an $F$-stable submodule, we have the desired equality. Putting everything together we see that $\Rlim_FM\simeq M^0[0]$. 

Thus taking cohomology in \eqref{isoinf} we get an isomorphism
\[
\Hinf{i}(X/\W(k))\simeq\bigcap_nF^n\left(\Hcris{i}(X/\W(k)\right)
\]
and the right hand side is the unit-root part of $\Hcris{i}(X/\W(k))$.
$\qed$

\section{Proof of Theorem \ref{mainthm}}\label{pruffa}


In this section we prove Theorem \ref{mainthmbl} following the argument given in \cite[Section 7]{bhattlurie}. We recall the statement in Theorem \ref{triangololiscio}. After defining the maps that appear in the statement there is a quite formal descent step, and the bulk of the work is then showing that the statement remains true for an eqrsp $\Fp$-algebra. We do this in Theorem \ref{triangolosemipft}. Once all of this is done Theorem \ref{triangololiscio} follows with just a little more effort. The first step is define a refinement of the crystalline first Chern class.

\begin{conssezione}
    
Let $X$ be an $\Fp$-scheme. There is a short exact sequence of sheaves on $\Cris(X/\Zp)$,
$$0\to(1+\icris)^{\times}\to\ocris^{\times}\to\Gm\to0,$$ 
which we view as a morphism $\Gm[-1]\to(1+\icris)^{\times}$ in the derived category. Note in passing that each local section of $\icris$ is locally nilpotent, so in fact $(1+\icris)^{\times}=1+\icris$. We will keep writing $(1+\icris)^{\times}$ to underline the multiplicative group structure.

The ideal $\icris$ has a canonical PD-structure, so we can define a logarithm 
\[
\log:(1+\icris)^{\times}\to\icris \quad\quad 1+x\mapsto\sum_{n\ge1}(-1)^{n+1}\delta_n(x)(n-1)!
\]
where $(\delta_n)$ is the PD-structure of the test object. By composing we get a map $\Gm[-1]\to\icris$ in the derived category, and taking derived global sections a map $c_1:\Rflat(X,\Gm)[-1]\to\Nyg\Rcris(X/\Zp)$.

By Lemma \ref{derivedpcompletion} and Example \ref{completamentoGm}, the map $c_1$ factors through the derived $p$-completion of $\Rflat(X,\Gm)[-1]$, inducing a map $\hat{c}_1:\Rflat(X,\Zp(1))\to\Nyg\Rcris(X/\Zp)$. This is the refinement of the crystalline first Chern class that we needed. The map comparing flat and crystalline cohomology in Theorem \ref{mainthm} is the composition 
\begin{equation}\label{comparisonmap}
\Rflat(X,\Zp(1))\xto{}{\hat{c}_1}\Nyg{}\Rcris(X/\W(k))\to\Rcris(X/\W(k)).
\end{equation}

\end{conssezione}

\begin{exsezione}\label{logqrsp}

Let us spell out this construction for $X=\Spec(C)$ with $C$ semiperfect, on the level of $H^0$. We have 
\begin{equation*}
   \Hflat{0}(X,\Zp(1))=\varprojlim_n\mupn(C)=T_p(C^\times),
\end{equation*}
and seeing $\mupn(C)$ as a subset of $C$ this group is naturally a subset of $C^\flat$. Specifically, it consists of the sequences $(\dots,a_1,a_0)\in\perf{C}$ such that $a_0=1$, $a_1^p=1$ \textit{etc.}, which is equivalently described as $1+\ker(C^\flat\onto{} C)\subseteq C^\flat$. On the other hand, $\Hcris{0}(X,\icris)=\Nyg{}\Acr(C)$.

The map $\Gm[-1]\to(1+\icris)^{\times}$ induces 
\[
\Gm[-1]\ten^L_{\Z}\Z/p^n=\mupn\to(1+\icris)^{\times}\ten^L_{\Z}\Z/p^n=(1+\icris)^{\times}/p^n(1+\icris)^{\times}.
\]
Note that $\cap_n(p^n(1+\icris)^\times)=\cap_n((1+\icris^{p^n})^\times)=1$, because all local sections of $\icris$ are nilpotent. Taking cohomology and the derived limit gives a map
\[
\Rflat(X,\Zp(1))\to\Rgam(X,(1+\icris)^{\times}),
\]
and $\hat{c}_1$ is the composition of this map with 
\[
\Rgam(X,(1+\icris)^{\times})\xto{}{\log}\Rgam(X,\icris).
\]

The map $\mupn\to(1+\icris)^{\times}/p^n(1+\icris)^{\times}$ comes from the snake lemma applied to the diagram
\[\begin{tikzcd}
	0 & {(1+\icris)^{\times}} & {\ocris^{\times}} & \Gm & 0 \\
	0 & {(1+\icris)^{\times}} & {\ocris^{\times}} & \Gm & 0
	\arrow[from=1-1, to=1-2]
	\arrow[from=1-2, to=1-3]
	\arrow["{p^n}", from=1-2, to=2-2]
	\arrow[from=1-3, to=1-4]
	\arrow["{p^n}", from=1-3, to=2-3]
	\arrow[from=1-4, to=1-5]
	\arrow["{p^n}", from=1-4, to=2-4]
	\arrow[from=2-1, to=2-2]
	\arrow[from=2-2, to=2-3]
	\arrow[from=2-3, to=2-4]
	\arrow[from=2-4, to=2-5]
\end{tikzcd}\]
On the level of global sections, $x\in\mupn(C^{\times})$ is sent to the class of $[y]^{p^n}$ where $y$ is any element of $\perf{C}$ lifting $x$. Taking the projective limit we get a map $T_p(C^\times)\to\Nyg{}\Acr(C)$ which maps $x$ to its Teichmüller representative $[x]\in\Acr(C)$. Therefore, $H^0(\hat{c}_1)$ maps $x\in T_p(C^{\times})$ to $\log([x])=-\sum_{n\ge1}\frac{([x]-1)^n}{n}\in\Acr(C)$.

\end{exsezione}

\begin{conssezione}\label{costruzionsfsup}

Let $A$ be an eqrsp $k$-algebra. By virtue of Corollary \ref{filn}, the map $\Nyg\Acr(A)\to\Acr(A)$ which maps $x$ to $F(x)/p-x$ is well-defined, so we obtain a map $F/p-1:\Nyg\Rcris(A/\Zp)\to\Rcris(A/\Zp)$. If $B$ is another eqrsp algebra with a map $f:A\to B$, the diagram 
\[\begin{tikzcd}
	{\Nyg\Rcris(A/\Zp)} & {\Rcris(A/\Zp)} \\
	{\Nyg\Rcris(B/\Zp)} & {\Rcris(B/\Zp)}
	\arrow["{F/p-1}", from=1-1, to=1-2]
	\arrow[from=1-1, to=2-1]
	\arrow[from=1-2, to=2-2]
	\arrow["{F/p-1}", from=2-1, to=2-2]
\end{tikzcd}\]commutes. This allows Bhatt-Lurie to define the map $F/p-1:\Nyg\Rcris(X/\Zp)\to\Rcris(X/\Zp)$ for any smooth $k$-scheme $X$ via descent, i.e. as the map $\Tot(F/p-1)$ between the diagrams of Corollary \ref{discesaesplicita}. However this procedure relies on the machinery of $\infty$-categories to guarantee that the map thus obtained is canonical, and does not depend on a chosen covering of $X$.

Instead, for a smooth $k$-scheme $X$ we will use the ideas of Section \ref{disc} to write down canonical and functorial complexes representing $\Rcris(X/\Zp)$ and $\Nyg{}\Rcris(X/\Zp)$, and we will define $F/p-1$ as a map between these two complexes. If $\mathfrak{U}=(U_i)$ is a finite affine cover of $X$, let $\Acr(\mathfrak{U}_{\mathrm{perf}})$ denote the naive \v{C}ech complex 
\[
\bigoplus_i\Acr(\coperf{(U_i)})\to\bigoplus_{i,j}\Acr(\coperf{(U_i)}\times_X\coperf{(U_j)})\to\dots
\]
According to the \v{C}ech spectral sequence and Corollary \ref{cocrsp}, this complex computes the crystalline cohomology of $X$. Denote by $\Acr(\coperf{X})$ the filtered colimit of the complexes $\Acr(\mathfrak{U}_{\mathrm{perf}})$, where there is a map for every cover refinement. The complex $\Acr(\coperf{X})$ is canonical, functorial in $X$ and computes $\Rcris(X/\Zp)$. 

Define likewise $\Nyg{}\Acr(\coperf{\mathfrak{U}})$ and $\Nyg{}\Acr(\coperf{X})$. Then there is a canonical map of complexes
\[
F/p-1:\Nyg{}\Acr(\coperf{X})\to\Acr(\coperf{X})
\]
obtained from the maps $F/p-1$ for quasiregular semiperfect algebras. This then defines a canonical map $F/p-1:\Nyg{}\Rcris(X/\Zp)\to\Rcris(X/\Zp)$, which is functorial for maps $f:Y\to X$ between smooth schemes, or between an eqrsp $Y$ and a smooth $X$. 
\end{conssezione}

We are now ready to state and prove the main theorems of the paper. For a $k$-scheme $X$ which is either smooth or semiperfect consider the triangle
\begin{equation}\label{triangolo}
    \Rflat(X,\Zp(1))\xto{}{\hat{c}_1}\Nyg\Rcris(X/\Zp)\xto{}{F/p-1}\Rcris(X/\Zp)
\end{equation}
in $\D(\Zp)$.

\begin{thmsezione}\label{triangolosemipft}

(\cite[Theorem 7.1.1]{bhattlurie}) Let $C$ be an eqrsp $\Fp$-algebra. The following hold:
\begin{enumerate}
    \item[(1)] For $X=\Spec(C)$, the triangle \eqref{triangolo} reduces to the sequence of abelian groups
    \begin{equation}\label{triangoloeqrsp}
    0\to(1+J)^{\times}\xto{}{\log}\Nyg\Acr(C)\xto{}{F/p-1}\Acr(C)\to0,
    \end{equation}
    where $J=\ker(\perf{C}\onto{}{}C)$. In other words, the natural map 
    $(1+J)^{\times}[0]\to\Rflat(X,\Zp(1))$ is a quasi-isomorphism in $\D(\Z)$.
    \item[(2)] The sequence \eqref{triangoloeqrsp} is exact.
\end{enumerate}

\end{thmsezione}

\noindent{\it Proof.} First we prove that (1) follows from (2). Assume that (2) holds. As was explained in Example \ref{logqrsp} we have $H^0(\Spec(C),\Zp(1))=\lim\mupn(C)=T_p(C^\times)$ which as a subset of $C^\flat$ is identified with $1+J$. Furthermore here $1+J$ and $(1+J)^\times$ are equal because $C^\flat$ is $J$-complete (Proposition \ref{ainfspqr}). Corollary \ref{cocrsp} gives $\Nyg{}\Rcris(X/\Zp)\simeq\Nyg{}\Acr(C)[0]$ and $\Rcris(X/\Zp)\simeq\Acr(C)[0]$, and under these identifications $H^0(\widehat{c}_1)$ becomes the logarithm map (see Example \ref{logqrsp}). Therefore if we prove that $\Hflat{i}(X,\Zp(1))=0$ for $i>0$ we have (1). To this end it is enough to show that $\Hflat{i}(X,\mu_p)=0$ for $i>0$.

The Kummer sequence $0\to\mup\to\Gm\xto{}{}\Gm\to0$ on the small étale site of $X$ is exact, because $C$ is semiperfect, so using the fact that étale and fppf cohomology of $\Gm$ coincide (see the proof of  Corollary \ref{discesagm}) we see that $\Hflat{i}(X,\mup)=\Het{i}(X,\mup)$ for all $i$. Thus we reduce to showing that $\Het{i}(X,\mup)=0$ for all $i>0$, and by Lemma \ref{discesaaciclica} it suffices to show that the functor $\mu_p(-)$ satisfies descent along étale covers of eqrsp algebras (recall from Proposition \ref{modelloetale} that an algebra which is étale over a eqrsp algebra is eqrsp). To show that $\mu_p$ satisfies descent along these étale covers it is enough to show the same for the functor $T_p(-^{\times})$. But now, assuming (2), the group $T_p(-^{\times})$ is the fibre of $\Nyg\Acr(-)\to\Acr(-)$. By Proposition \ref{discesaetalecris} these two groups satisfy descent along étale covers of eqrsp algebras, so $T_p(-^{\times})$ does as well. This proves (1).

Let us prove (2), i.e. that \eqref{triangoloeqrsp} is exact, in the case where $C=B[x^{\pmeninf}]/(x)$ for some perfect algebra $B$. The proof adapts as is to the general case, it is only more tedious to keep track of all indices. All groups in the sequence are $p$-torsion-free and $p$-complete, so it is enough to show that \eqref{triangoloeqrsp} is exact mod $p$. Explicitly, we have a sequence
\begin{equation}\label{triangolomodp}
    0\to(1+J)^\times/p\xto{}{\overline{\log}}\Nyg\Acr(C)/p\xto{}{M}\Acr(C)/p\to0
\end{equation}
and the following descriptions of these groups.
\begin{enumerate}
    \item[(i)] Since $C^\flat$ is $J$-complete (Proposition \ref{ainfspqr}), we have $(1+J)^\times=(1+J)$, so that $p(1+J)^\times=(1+J)^p=(1+J^p)$. Using the identification $(1+J)^\times=T_p(C^\times)$, we see that any element of $(1+J)^{\times}/(1+J^p)=\mu_p(C^\times)$ can be written uniquely as $1+\sum b_{\alpha}x^{\alpha}$, where the sum is finite and runs over the $\alpha$ such that $1\le\alpha<p$.
    \item[(ii)] $\Nyg\Acr(C)/p$ is the free $B$-module with basis
    $$\{px^{\alpha}\text{ s.t. }\lfloor\alpha\rfloor=0\}\cup\left\{\frac{x^{\alpha}}{(\alpha!)_p}\text{ s.t. }\lfloor\alpha\rfloor\ge1\right\}$$
    \item[(iii)] $\Acr(C)/p$ is the free $B$-module with basis all $\frac{x^{\alpha}}{(\alpha!)_p}$
    \item[(iv)] For all $x\in J$:
    $$\log([1-x])=\sum_{d=1}^p(-1)^{d-1}\frac{[x]^d}{d}$$
    in $\Acr(C)/p$.
    \item[(v)] For all $b\in B$ and all $\alpha$:
    $$M(bpx^{\alpha})=b^px^{p\alpha}\hspace{3.05cm}\text{if }\lfloor\alpha\rfloor=0$$
    $$M\left(b\frac{x^{\alpha}}{(\alpha!)_p}\right)=\begin{cases}
        b^p\frac{x^{p\alpha}}{(p\alpha!)}_p-b\frac{x^{\alpha}}{(\alpha!)}_p&\text{if }{\lfloor\alpha\rfloor}=1\\
        -b\frac{x^{\alpha}}{(\alpha!)_p}&\text{if }{\lfloor\alpha\rfloor\ge2}
    \end{cases}$$
\end{enumerate}
Points (i), (ii), (iii) and (v) are easily derived from the explicit descriptions of $\perf{C},\Acr(C)$ and $\Nyg\Acr(C)$ given in Subsection \ref{elqrsp}. Point (iv) is true because the equality $[a+b]=[a]+[b]$ holds mod $p$.

\textbf{Exactness on the left:} take $a=1+\sum b_{\alpha}x^{\alpha}$ in $(1+J)^{\times}/(1+J^p)$. Suppose $a\ne1$ and let $\alpha_0$ be the smallest $\alpha$ such that $b_{\alpha}\ne0$. Then the coefficient in front of $\frac{x^{\alpha_0}}{(\alpha_0!)_p}$ for $\overline{\log}(a)$ is $b_{\alpha_0}$, which is nonzero by hypothesis. Therefore $a\notin\ker(\overline{\log})$ so $\overline{\log}$ is injective.

\textbf{Exactness in the middle:} let $m\in \ker M$, and write
\begin{equation*}
    m=\sum_{\alpha<1}b_\alpha px^\alpha+\sum_{\alpha\ge1}b_\alpha\frac{x^\alpha}{(\alpha!)_p}
\end{equation*}
in the basis of (ii), with the $b_\alpha$ in $B$. Let $\alpha_m$ be the largest $\alpha$ such that $b_\alpha\ne0$. If $\alpha_m>2p$, the coefficient of $\frac{x^{\alpha_m}}{(\alpha_m!)_p}$ in $M(m)$ is $b_{\alpha_m}$, which is non-zero. This contradicts $M(m)=0$, so $\alpha_m<2p$. Similarly, if $\alpha_m<2$ the coefficient of $\frac{x^{p\alpha_m}}{((p\alpha_m)!)_p}$ in $M(m)$ is $b^p_{\alpha_m}$, which is non-zero. This contradicts $M(m)=0$. In conclusion, $2\le \alpha_m<2p$.

This discussion implies that it suffices to check the surjectivity of the composition
$$
\theta:(1+J)^{\times}/(1+J^p)\to\Nyg\Acr(C)/p\onto[]{}S,
$$
which maps the class of $1+x\in (1+J)^\times$ to $\log([1+x])$, and then projects onto $S$, the free $B$-module with basis $\left\{\frac{x^{\alpha}}{(\alpha!)_p}\text{ s.t. }2\le\alpha<2p\right\}$, thus forgetting all coefficients $b_\alpha$ of $\log([1-x])$ with $\alpha<p$ and $\alpha\ge2p$. Write furthermore $S=S_1\oplus S_2$, where $S_1$ (resp. $S_2$) is the free $B$-module with basis the $\frac{x^\alpha}{(\alpha!)_p}$ such that $2\le\alpha<p$ (resp. $p\le\alpha<2p$). In particular $S_1=0$ for $p=2$.

Take $a=1+\sum b_{\alpha}x^{\alpha}$ in $(1+J)^{\times}/(1+J^p)$. For $p\le\alpha<2p$ the coefficient in front of $\frac{x^{\alpha}}{(\alpha!)_p}=\frac{x^{\alpha}}{p}$ for $\overline{\log}(a)\in\Nyg{}\Acr(C)/p$ is equal to $b_{\alpha/p}^p$. Thus the composition of $\theta$ with the projection $S\onto{}{}S_2$ is surjective, and if $p=2$ we are done. In the case $p\ne2$, note that if $b_{\alpha}=0$ for $1\le\alpha<2$, which is equivalent to $a\in(1+J^2)^{\times}/(1+J^p)$, then the projection of $\theta(a)$ onto $S_2$ is $0$. We may thus reduce to showing that the composition $$(1+J^2)^{\times}/(1+J^p)\to(1+J)^{\times}/(1+J^p)\to S\onto[]{}S_1$$ is surjective. 

The point now is that there is a natural identification $J^2/J^p\simeq S_1$ - where we view $J^2$ as a $B$-submodule of $\Acr(C)/p$ - and the resulting map $(1+J^2)^{\times}/(1+J^p)\to J^2/J^p$ sends $1+x$ to $\sum_{d=1}^{p-1}(-1)^d\frac{x^d}{d}$. An explicit inverse of this map is given by $x\mapsto\sum_{d=1}^{p-1}\frac{x^d}{d!}$.

\textbf{Exactness on the right:} let $b\in B$, $\alpha\in\Z[1/p]_+$, $a=b\frac{x^{\alpha}}{(\alpha!)_p}$, and let's prove that $a\in\im(M)$. If $\lfloor\alpha\rfloor\ge2$ then $M(-a)=a$ so $a\in\im(M)$. If $\lfloor\alpha\rfloor=0$ then $M(b^{1/p}px^{\alpha/p})=a$ so $a\in\im(M)$. Finally, if $\lfloor\alpha\rfloor=1$ then $M(a)+a=b^p\frac{x^{p\alpha}}{((p\alpha)!)_p}$ and the latter lies in $\im(M)$ because $\lfloor p\alpha\rfloor\ge2$. Thus in this case too $a\in\im(M)$. $\qed$

\begin{thmsezione}\label{triangololiscio}

(\cite[Theorem 7.3.5]{bhattlurie}) If $X$ is smooth over a perfect field $k$ then triangle \eqref{triangolo} is exact.

\end{thmsezione}

\noindent{\it Proof.} If all functors in \eqref{triangolo} satisfy descent along $Y\to X$, and $\eqref{triangolo}$ is exact for $Y,Y\times_XY,$ etc. then it is exact for $X$. By Zariski descent, Proposition \ref{discesaetalecris} and Corollary \ref{discesagm}, we thus reduce to an affine $X$. By descent along $\coperf{X}\to X$, Theorem \ref{discesafrobsmooth} and Lemma \ref{discesaperfgm}, we reduce to an eqrsp scheme. Then the result is true by Theorem \ref{triangolosemipft}.
$\qed$

\noindent{\it Proof of Theorem \ref{mainthm}.}
Now $k$ is algebraically closed, and we revert to writing $\Hcris{i}(X/\W(k))$ instead of $\Hcris{i}(X/\Zp)$ to highlight that we are dealing with $\W(k)$-modules. Let $X$ be smooth and proper. Taking cohomology of \eqref{triangolo} we get the exact sequence
\begin{align*}
    \Nyg&\Hcris{i-1}(X/\W(k))\xto{}{F/p-1}\Hcris{i-1}(X/\W(k))\to \Hflat{i}(X,\Zp(1))\to\dotsb\\
    &\dotsb\to\Nyg\Hcris{i}(X/\W(k))\xto{}{F/p-1}\Hcris{i}(X/\W(k))
\end{align*}
We take as our comparison map the composition 
\begin{equation*}
\alpha_X:\Hflat{i}(X,\Zp(1))\xto{}{\hat{c}_1}\Nyg\Hcris{i}(X/W(k))\to\Hcris{i}(X/\W(k))
\end{equation*}
as in \eqref{comparisonmap}. In Example \ref{defnygaard} we saw that the natural maps 
$\Nyg\Hcris{i}(X/\W(k))\to\Hcris{i}(X/\W(k))$ become isomorphisms after inverting $p$, so we get the long exact sequence
\[
\dots\xto{}{F/p-1}\Hcris{i-1}(X/\W(k))[1/p]\to \Hflat{i}(X,\Qp(1))\xto{}{\alpha_X}\Hcris{i}(X/\W(k))[1/p]\xto{}{F/p-1}\dots
\]
Since $X$ is proper, the crystalline cohomology groups are finite-dimensional $\W(k)[1/p]$-vector spaces. It then follows from the Dieudonné-Manin classification that each of the maps $\Hcris{j}(X/\W(k))[1/p]\xto{}{F/p-1}\Hcris{j}(X/\W(k))[1/p]$ is surjective, see \cite[5.6 in II.5.A]{illusiedrw}. Thus the sequence $$0\to \Hflat{i}(X,\Qp(1))\xto{}{\alpha_X}\Hcris{i}(X/\W(k))[1/p]\xto{}{F/p-1}\Hcris{i}(X/\W(k))[1/p]\to0$$ is exact, and $\alpha_X$ gives the desired isomorphism 
\[
\Hflat{i}(X,\Qp(1))\simeq\left(\Hcris{i}(X/\W(k))[1/p]\right)^{F=p}.
\]
$\qed$

\begin{remsezione}\label{comparisonetale}

The same arguments apply to show that if $X$ is a smooth scheme over a perfect field $k$, there is an exact triangle
\[
\Ret(X,\Zp)\to\Rcris(X/\Zp)\xto{}{F-1}\Rcris(X/\Zp)
\]
where $\Ret(X,\Zp)=\Rlim_n\Ret(X,\Z/p^n)$. By descent one reduces to showing that if $C$ is eqrsp and $\Spec(C)$ is connected, the sequence
\[
0\to\Zp\to\Acr(C)\xto{}{F-1}\Acr(C)\to0
\]
is exact, which is much easier than Theorem \ref{triangolosemipft}. It follows that when $X$ is proper over the perfect field $k$ there is an isomorphism
\[
\Het{i}(X,\Zp)\simeq\Hcris{i}(X/\W(k))^{F=1}
\]
and $\Het{i}(X,\Zp)\otimes_{\Zp}\W(k)$ is isomorphic to the unit-root part of $\Hcris{i}(X/\W(k))$

\end{remsezione}

\begin{remsezione}\label{}

Theorem \ref{discesafrobsmooth} and Theorem \ref{triangolosemipft}, are proved in greater generality in \cite{bhattlurie}. Recall that a morphism of $\Fp$-algebras $f:A\to B$ is quasisyntomic if it is flat and $ \Ll_{B/A}$ has Tor-amplitude in $[-1,0]$. If $f$ is faithfully flat it is a quasysintomic cover, and an $\Fp$-algebra $A$ is quasisyntomic if $\Fp\to A$ is a quasisyntomic morphism. Frobenius-smooth algebras and eqrsp algebras are quasisyntomic, and if $A$ is Frobenius-smooth then the map $A\to\coperf{A}$ is a quasisyntomic cover. 

Bhatt--Lurie prove that $\Rcris(-/\Zp)$ satisfies descent along quasisyntomic covers of quasisyntomic algebras \cite[Corollary 4.6.7]{bhattlurie}. A quasisyntomic $\Fp$-algebra $A$ has a cover $A\to\coperf{A}$ such that all products $B,B\otimes_A B,\cdots$ are qrsp (Definition \ref{spqrdef}), and they prove that triangle \eqref{triangolo} is exact for qrsp algebras.

These results are more general and more natural than the ones presented here, at the expense of being less elementary. The conceptual advantage in their approach is to consider the quasisyntomic topology as a Grothendieck topology for which crystalline cohomology satisfies descent, which is more satisfying than our ad hoc approach. Their formulation of descent in the derived category is made possible by working with the derived $\infty$-categories.

\end{remsezione}


\section{Applications to fppf cohomology}\label{applicazioni}

In this final section we give some applications of Theorem \ref{triangololiscio} to the structure and properties of the groups $\Hflat{i}(X,\Zp(1))$, which, as we will see, behave quite differently from their étale $\ell$-adic counterparts $\Het{i}(X,\Z_\ell(1))$. In Subsection \ref{sezionefinitezza} we show that if $X$ is smooth and proper over $k$, and $k$ is algebraically closed, then the groups $\Hflat{i}(X,\Zp(1))$ are the direct sum of a finite free $\Zp$-module, and of a $p$-torsion group of finite $p$-exponent. This $p$-torsion group may be infinite, even for abelian varieties, but for $i=1,2$ it is always of finite type. In Subsection \ref{straight} we show that when $X$ is sufficiently nice (e.g. an abelian variety) fppf cohomology is uniquely determined by crystalline cohomology. We use this to compute the action of the multiplication-by-$n$ map $[n]$ on the fppf cohomology of an abelian variety, a problem suggested to the author by A.Skorobogatov. In the end we show via two examples how the Nygaard filtration can be used to do explicit computations. 

Many of the results presented here are well-known: except for the action of $[n]$ on fppf cohomology, they appear in papers by Illusie \cite{illusiedrw}, Illusie-Raynaud \cite{illusieraynaud} and Milne \cite{milnezeta}, sometimes with more precise statements. But as far as the author knows the proofs presented here are new in the literature.

For the reader's convenience, we recall here the main properties of the Nygaard filtration and its relation to fppf cohomology. We suppose in the whole section that $k$ is algebraically closed, and that $X$ is a proper smooth variety over $k$. We will write $\Rcris(X)$ and $\Hcris{i}(X)$ for crystalline cohomology, thus omitting the coefficient ring $\W(k)$, as we believe this will not cause any confusion.

By definition, the Nygaard filtration sits in an exact triangle
\begin{equation}\label{triangoloagain}
\Nyg{}\Rcris(X)\to\Rcris(X)\to\Rgam(X,\calO_X)
\end{equation}
in $\D(\Zp)$. In other words there is a long exact sequence
\begin{equation}\label{sel2}
	\cdots \to {H^{i-1}(X,\calO_X)} \to {\Nyg\Hcris{i}(X)} \to {\Hcris{i}(X)} \to {H^i(X,\calO_X)} \to \cdots,
\end{equation}
where ${\Hcris{i}(X)} \to {H^i(X,\calO_X)}$ is the canonical projection, which factors as 
\[
\Hcris{i}(X)\to\HdR{i}(X/k)\to H^{i}(X,\calO_X).
\]
It follows from the definition that $\Nyg{}\Hcris{i}(X)[1/p]\simeq\Hcris{i}(X)[1/p]$. 

The main result Theorem \ref{triangololiscio} asserts that there is a long exact sequence:
\begin{equation}\label{sel1}
	\cdots\to{\Hflat{i}(X,\Zp(1))}\to{\Nyg\Hcris{i}(X)}\xto{}{F/p-1}{\Hcris{i}(X)}\to{\Hflat{i+1}(X,\Zp(1))}\to\cdots
\end{equation}
The maps $F/p-1$ are not $\W(k)$-linear, they are the difference of a $\sigma$-linear map $F/p$ and a linear map $1$, which is the canonical map $\Nyg\Hcris{i}(X)\to\Hcris{i}(X)$. By a standard fact in the theory of isocrystals, \cite[II.5.A.5.3]{illusiedrw}, the maps $F/p-1$ are surjective after inverting $p$, so in particular $\coker(F/p-1)$ is a $p$-torsion group.

\subsection{Finiteness of fppf cohomology}\label{sezionefinitezza}

Our first goal is to prove the following result.

\begin{thm}\label{finiteexponent}

For all $i\ge0$ the $\Zp$-module $\Hflat{i}(X,\Zp(1))$ is isomorphic to the direct sum of $\Hflat{i}(X,\Zp(1))_{\text{tors}}$ with a free $\Zp$-module of rank $r=\rk_{\W(k)}\left(\Hcris{i}(X/\W(k))^{F=p}\right)$. Moreover, $\Hflat{i}(X,\Zp(1))_{\text{tors}}$ is a $p$-group of finite $p$-exponent.

\end{thm}

We break up the proof into smaller statements.

\begin{lem}\label{lemma1}
    Let $i\ge0$. The torsion of $\ker\left(F/p-1:\Nyg\Hcris{i}(X)\to\Hcris{i}(X)\right)$ is a $p$-group of finite $p$-exponent.
\end{lem}

\noindent{\it Proof.} The torsion of $\ker\left(F/p-1:\Nyg\Hcris{i}(X)\to\Hcris{i}(X)\right)$ is a subgroup of the torsion subgroup of $\Nyg\Hcris{i}(X)$, which is a finite type $\W(k)$-module, and thus has finite $p$-exponent. This proves the first statement.$\qed$

\begin{lem}\label{lemma2}
    Let $i\ge0$. The group $\coker\left(F/p-1:\Nyg\Hcris{i}(X)\to\Hcris{i}(X)\right)$
    is a $p$-group of finite $p$-exponent. 
\end{lem}

\noindent{\it Proof.} We reduce the problem to one of semilinear algebra. The torsion of $\Nyg\Hcris{i}(X)$ and of $\Hcris{i}(X)$ is of finite $p$-exponent, so we may work modulo torsion. Taking $M=\Hcris{i}(X)/\text{tors}$ and $N=\Nyg{}\Hcris{i}(X)/\text{tors}$, we are in the following situation: let $M$ be a finite free $\W(k)$-module with a $\sigma$-linear endomorphism $F$, and $N\subseteq M$ an $F$-stable submodule containing $pM$ such that $F(N)\subseteq pM$. We want to prove that the cokernel of $F/p-1:N\to M$ is a torsion group of finite $p$-exponent. It is enough to prove that $F/p-1:pM\to M$ has cokernel of finite $p$-exponent, or that $F-p:M\to M$ has cokernel of finite $p$-exponent. Upon replacing $M$ with a submodule of the same rank, we may suppose by the Dieudonné-Manin theorem that $M$ is a direct sum of modules $M_{s/r}=\W(k)^r$ where $F(x_1,\dots,x_r)=(x_2,\dots,x_{r},p^sx_1)$, and further reduce to a single such $M_{s/r}$ (here $r>0,s\ge0$ are coprime integers). Now we distinguish three cases.
\begin{itemize}
    \item[(1)]If $r=s=1$ then $F-p$ is surjective onto $pM$: $(F-p)(\lambda x_1)=x_1\left(p\sigma(\lambda)-p\lambda\right)$ and $\lambda\mapsto p\sigma(\lambda)-p\lambda$ is surjective onto $p\W(k)$.
    \item[(2)]If $s/r>1$ then $F-p$ is surjective onto $p^rM$: for $a\in p^rM$ it is easy to see that $F^i(a)/p^i$, which is an element of $M[1/p]$, actually belongs to $M$ for all $i\ge1$, that the infinite sum $a+\sum_{i\ge1}F^i(a)/p^i$ converges to some $x\in M$ and that $F(x)-px=a$.
    \item[(3)]If $s/r<1$ then $F-p$ is surjective onto $p^rM$: once again for $a\in p^sM$ one checks that the infinite sum $\sum_{i\ge1}F^{-i}(p^{i-1}a)$ makes sense in $M$ and converges to some element $x\in M$ such that $F(x)-px=a$.$\qed$
\end{itemize}

We need one last technical lemma before tackling Theorem \ref{finiteexponent}.

\begin{lem}\label{lemmapdiv}

Let $L_0\to L_1\to L_2\to\dots$ be a sequence of $\W(k)$-modules. If the $L_i$ are free of the same rank and the transition maps are injective but not surjective, the colimit $L=\varinjlim_nL_n$ contains a non-trivial $K$-vector space.

\end{lem}

\noindent{\it Proof.} We view the $L_i$ as strictly increasing sub-$\W(k)$-modules of the $K$-vector space $L_0[1/p]$, and $L$ as their union. Take a sequence $x_n$ of elements of $L_0$ such that each $x_n$ is not divisible by $p$ in $L_0$ but is divisible by $p^n$ in $L$. Such $x_n$ exist because the length of $L_i/L_0$ increases strictly with $i$. Let 
\[
E=\bigcap_m\text{Span}_K(\{x_n,n\ge m\}),
\]
i.e. for any $m>0$ there is a basis of $E$ made from a subset of $\{x_n,n\ge m\}$. This vector space is non-zero, and let $E_0=E\cap L_0$. Then for every $m>0$ the $\W(k)$-module $L$ contains $\{x_n/p^m,n\ge m\}$, so $L$ contains $1/p^mE_0$ and also their union which is the whole of $E$.$\qed$

\noindent{\it Proof of Theorem 5.1.1}
For $i=0$ the statement is immediate. Fix $i\ge1$, let $r=\rk_{\W(k)}\left(\Hcris{i}(X)^{F=p}\right)$ and define $M$ by the short exact sequence
\[
0\to \Hflat{i}(X,\Zp(1))_{\text{tors}}\to\Hflat{i}(X,\Zp(1))\to M\to0,
\]
so that $M$ is a torsion-free $\Zp$-module. By Theorem \ref{mainthm} $M[1/p]$ is isomorphic to $\Qp^r$, so if $M$ were not a finite type $\Zp$-module it would contain an infinite strictly increasing chain of finite free $\Zp$-submodules $L_0\subset L_1\subset\cdots\subset M$.  Their rank is bounded above by $r$, so by Lemma \ref{lemmapdiv} $M$ must contain an element infinitely divisible by $p$. But $M$ injects into the free $\W(k)$-module $\Nyg\Hcris{i}(X)/\text{tors}$, so this yields a contradiction. Thus $M\simeq\Zp^r$ and the first part of the statement follows.

It remains to prove that $M_{\text{tors}}$ is a $p$-group of finite $p$-exponent. It lies in a short exact sequence coming from \eqref{sel2}:
\[
0\to\coker\left(F/p-1\right)\to M_{\text{tors}}\to\ker\left(F/p-1\right)_{\text{tors}}\to0
\]
where the first map is $F/p-1$ in degree $i-1$ and the second map is $F/p-1$ in degree $i$ - note that we have used here that the cokernel of $F/p-1$ is torsion. To conclude, it is enough to show that $\ker\left(F/p-1\right)_{\text{tors}}$ and $\coker\left(F/p-1\right)$ are $p$-groups of finite $p$-exponent. This is the content of Lemma \ref{lemma1} and Lemma \ref{lemma2} respectively.$\qed$

As mentioned in the introduction, this result can be applied to the structure of to the $p$-torsion of the Brauer group.

\begin{cor}\label{finitebrauer}

The group $\Br(X)[p^{\infty}]$ is the direct sum of a divisible $p$-group of 
the form $\left(\Qp/\Zp\right)^a$, where $a=\dim_{\Qp}\Hflat{2}(X,\Qp(1))-\rk\NS(X)$, 
and of a $p$-group of finite $p$-exponent.

\end{cor}

\noindent{\it Proof.} From the short exact sequence \eqref{brauer} mentioned in the introduction we see that $\Br(X)[p^\infty]$ is isomorphic to the direct sum of $(\Qp/\Zp)^a$ with $\Hflat{3}(X,\Zp(1))_\text{tors}$. The latter has finite $p$-exponent by Theorem \ref{finiteexponent}.
$\qed$

\begin{rem}\label{unipotente}

Theorem \ref{finiteexponent} admits the following generalization: for all $i\ge0$ the cokernel of $F/p-1:\Nyg{}\Hcris{i}(X)\to\Hcris{i}(X)$ is of the form $U_i(k)$, where $U_i$ is a connected unipotent algebraic $k$-group. Combined with Lemma \ref{lemma1}, this gives a short exact sequence
\[
0\to U_i(k)\to\Hflat{i}(X,\Zp(1))_{\text{tors}}\to F_i\to0
\]
with $F_i$ a finite group. The dimension of $U_i$ is easily bounded by the sum of Hodge numbers $h^{0,i-1}+h^{0,i}$, see \cite[p.554]{artinssk3}, and can be computed in terms of invariants of the de Rham-Witt complex. A statement encompassing these results can be found in \cite[IV.3.3]{illusieraynaud} - notice that $\Hflat{i}(X,\Zp(1))\simeq H^{i+1}(X,W\Omega^1_{\log})$, \cite[p.192]{illusieraynaud}. This implies in particular that if $U_i$ is nonzero then $\Hflat{i}(X,\Zp(1))_{\text{tors}}$ is necessarily infinite. In Subsection \ref{straight} we show that if $X=E\times E$ where $E$ is a supersingular elliptic curve then $U_3\simeq\Ga$.
\end{rem}

The rest of this section is devoted to showing that the groups $\Hflat{i}(X,\Zp(1))$ are finitely generated $\Zp$-modules when $i=1,2$. In order to use the sequence \eqref{sel1} effectively it is crucial to understand crystalline cohomology and the Nygaard filtration in degrees $0,1$. For degree $1$, Proposition \ref{comparisonalbanese} tells us that $\Hcris{1}(X)$ has the structure of a Dieudonné module, i.e. a free $\W(k)$-module with a $\sigma$-linear endomorphism $F$ and a $\sigma^{-1}$-linear endomorphism $V$ such that $FV=VF=p$. For degree $0$ we have the following lemma.

\begin{lem}\label{lem0}
    The maps $\Hcris{0}(X)/p\to\HdR{0}(X/k)\to H^0(X,\calO_X)$ are isomorphisms, and thus $\Nyg{}\Hcris{0}(X)=p\Hcris{0}(X)$.
\end{lem}

\noindent{\it Proof.}
By Proposition \ref{confrontoderham} there is a short exact sequence
\begin{equation*}
    0\to\Hcris{0}(X)/p\to\HdR{0}(X/k)\to\Hcris{1}(X)[p]\to0,
\end{equation*}
and we know that $\HdR{0}(X/k)=H^0(X,\calO_X)=k$. Choosing a $k$-point of $X$, we get $\Hcris{0}(\Spec(k))=\W(k)$ as a direct summand of $\Hcris{0}(X)$, and thus by comparing ranks with de Rham cohomology we find $\Hcris{0}(X)\simeq\W(k)$. The assertions of the lemma now follow easily
$\qed$
    
\begin{lem}\label{finitedegreeone}

The group $\Hflat{1}(X,\Zp(1))$ is a free $\Zp$-module of finite type.

\end{lem}

\noindent{\it Proof.}
By Theorem \ref{finiteexponent} it suffices to show that $\Hflat{1}(X,\Zp(1))$ is torsion-free. The map $\sigma/p-1:p\W(k)\to\W(k)$ is surjective. Therefore, combining Lemma \ref{lem0} with the sequence \eqref{sel1} we see that $\Hflat{1}(X,\Zp(1))$ is a subgroup of $\Nyg{}\Hcris{1}(X)$. By Lemma \ref{lem0} $\Hcris{0}(X)\to H^0(X,\calO_X)$ is surjective, thus \eqref{sel2} tells us  that the group $\Nyg{}\Hcris{1}(X)$ is a subgroup of $\Hcris{1}(X)$ . This last group is torsion-free, so we are done.
$\qed$

The case $i=2$ requires more work, and we first recall some basic properties of de Rham cohomology in positive characteristic. By functoriality, the absolute Frobenius of $X$ induces a map
\[
F:\HdR{i}(X/k)\to\HdR{i}(X/k)
\]
which is semi-linear with respect to $k$. The first observation is that this $F$ factors through the augmentation $\HdR{i}(X/k)\to H^{i}(X,\calO_X)$, thus defining a map $\overline{F}:H^{i}(X,\calO_X)\to\HdR{i}(X/k)$. Indeed, the pullback of a differential form by the absolute Frobenius is zero, so that the Frobenius map $F:\Omega_X^\bullet\to\Omega_X^\bullet$ factors through the map of complexes
\[
\calO_X[0]\to\Omega^{\bullet}_{X}
\]
given by the $p$-th power map in degree $0$. This last description makes it clear that $\overline{F}:H^0(X,\calO_X)\to\HdR{0}(X/k)$ is a ring homomorphism (the absolute Frobenius).

As we will now see, this map can also be described as a differential of the conjugate spectral sequence, via the Cartier isomorphism. Recall that if we equip the de Rham complex $\Omega_X^\bullet$ with the canonical filtration $\tau_{\le p}\Omega_X^\bullet$, we get a convergent spectral sequence
\begin{equation*}
    E^{p,q}_2=H^p(X,\mathcal{H}^q(\Omega^\bullet_X))\implies\HdR{p+q}(X/k)
\end{equation*}
called the conjugate spectral sequence. The conjugate spectral sequence is complemented by the Cartier isomorphism (see \cite[0.2.1.9]{illusiedrw}), which is an isomorphism of sheaves on $X^{(p)}$
\begin{equation*} 
    C^{-1}: \Omega_{X^{(p)}/k}^q \xrightarrow{\sim} \mathcal{H}^q(F_*\Omega_{X/k}^\bullet),
\end{equation*}
that induces semilinear isomorphisms $H^p(X,\mathcal{H}^q(\Omega^\bullet_X))\simeq H^ p(X,\Omega^q_X)$. In particular, we see that the Hodge-de Rham spectral sequence degenerates if and only if the conjugate spectral sequence degenerates, because both conditions are equivalent to $\sum_{p+q=n}h^{p,q}(X)=\dim\HdR{n}(X/k)$.

\begin{lem}\label{bordoconiugata}

Let $E_2^{\bullet,\bullet}$ denote the second page of the conjugate spectral sequence converging to $\HdR{i}(X/k)$. The Cartier isomorphism provides a semi-linear isomorphism $H^i(X,\calO_X)\simeq E_2^{i,0}$, and the composition
\begin{equation}\label{ultimaeq}
H^i(X,\calO_X)\simeq E_2^{i,0}\to\HdR{i}(X/k)
\end{equation}
coincides with $\overline{F}$
\end{lem}

\noindent{\it Proof.} We have $E^{i,0}_2=H^i(X,\mathcal{H}^0)$, where $\mathcal{H}^0$ is the kernel of $d:\calO_X\to\Omega^1_X$. By construction of the conjugate spectral sequence, the map $E^{i,0}_2\to\HdR{i}(X/k)$ is induced by the map of complexes $\mathcal{H}^0[0]\to\Omega^{\bullet}_{X}$. The Cartier isomorphism in degree $0$ says that $\calO_X\xto{}{F}\mathcal{H}^0$ is an isomorphism of sheaves, so the statement follows from the second description of $\overline{F}$ given above. $\qed$
\begin{lem}\label{comportamentofbar}

The following hold:
\begin{enumerate}
    \item The map $\overline{F}:H^{1}(X,\calO_X)\to\HdR{1}(X/k)$ is injective,
    \item If $\HdR{1}(X/k)\to H^1(X,\calO_X)$ is surjective and every global differential $1$-form is closed, then $\overline{F}:H^{2}(X,\calO_X)\to\HdR{2}(X/k)$ is injective.
\end{enumerate}
\end{lem}

\noindent{\it Proof.}
    Point (1) is true because trivially $E_2^{1,0}=E_{\infty}^{1,0}$. 
    
    In (2) we are assuming that 
    \[
    \dim_k\HdR{1}(X/k)=\dim_k H^{0}(X,\Omega^1_X)+\dim_k H^{1}(X,\calO_X),
    \]
    so by counting dimensions we find that have $E^{0,1}_2=E^{0,1}_{\infty}$ and $E^{1,0}_2=E^{1,0}_{\infty}$. Therefore the only differential involving $E^{2,0}_2$ is zero, i.e. $E^{2,0}_2=E^{2,0}_{\infty}$, which is what we needed to prove.
$\qed$

Let's look more closely at $\HdR{1}(X/k)$. Recall from Proposition \ref{crisabeliana} and Proposition \ref{comparisonalbanese} that there is a $\sigma^{-1}$-linear endomorphism $V$ of $\Hcris{1}(X)$ satisfying $F\circ V=V\circ F=p$.

\begin{lem}\label{diagrammacrisderhamlem}

There is a commutative diagram
\begin{equation}\label{diagrammacrisderham}
    \begin{tikzcd}
	{\Hcris{1}(X)} & {\Hcris{1}(X)} & {\Hcris{1}(X)} \\
	{H^1(X,\calO_X)} & {\HdR{1}(X/k)} & {\HdR{1}(X/k)}
	\arrow["F", from=1-1, to=1-2]
	\arrow["{q_1}"', from=1-1, to=2-1]
	\arrow["V", from=1-2, to=1-3]
	\arrow["q"', from=1-2, to=2-2]
	\arrow[from=1-3, to=2-3]
	\arrow["{\overline{F}}", hook, from=2-1, to=2-2]
	\arrow["V", from=2-2, to=2-3]
\end{tikzcd}
\end{equation}
where $V:\HdR{1}(X/k)\to\HdR{1}(X/k)$ is $\sigma^{-1}$-linear and factors through $H^0(X,\Omega_X^1)_{d=0}$.
\end{lem}
\noindent{\it Proof.}
The statements regarding the bottom row are proved in \cite[Proposition 5.7]{oda}. What is left to prove is the commutativity of the right-hand square. By Proposition \ref{comparisonalbanese} one reduces to the case of an abelian variety $A$. Then \cite[Corollary 5.11]{oda} shows that $\HdR{1}(A/k)$ is naturally identified with the Dieudonné module of $A[p]$ so the claim follows from Proposition \ref{crisabeliana}.
$\qed$

\begin{lem}\label{nygverschiebung}
We have the equality
\[
\Nyg{}\Hcris{1}(X)=V\left(\Hcris{1}(X)\right).
\]
\end{lem}

\noindent{\it Proof.}
The composition $\HdR{1}(X)\xto{}{V}\HdR{1}(X)\to H^1(X,\calO_X)$ is zero, so by diagram \eqref{diagrammacrisderham} we have the inclusion
\[
V\left(\Hcris{1}(X)\right)\subseteq\ker\left(\Hcris{1}(X)\to H^{1}(X,\calO_X)\right)=\Nyg{}\Hcris{1}(X),
\]
where the second equality follows from Lemma \ref{lem0} and the sequence \eqref{sel2}. Let us prove the other inclusion: if $a\in\Nyg{}\Hcris{1}(X)$, the existence of the map $F/p$ implies that $p$ divides $F(a)$. Thus we can write $F(a)=pb$ with $b\in\Hcris{1}(X)$, and applying $V$ we find $pa=pV(b)$. Since $\Hcris{1}(X)$ is torsion-free we have $a=V(b)$, which is what we wanted to prove.
$\qed$

\begin{lem}\label{fsupcoerente}

For all $i\ge0$, the composition
\begin{equation}\label{comp1}
H^i(X,\calO_X)\to\Nyg{}\Hcris{i+1}(X)\xto{}{F/p-1}\Hcris{i+1}(X)
\end{equation}
coming from the sequences \eqref{sel1} and \eqref{sel2}, and the composition
\begin{equation}\label{comp2}
H^i(X,\calO_X)\xto{}{\overline{F}}\HdR{i}(X/k)\to\Hcris{i+1}(X)[p]\subseteq\Hcris{i+1}(X)
\end{equation}
coming from Proposition \ref{confrontoderham}, are equal.
\end{lem}

\noindent{\it Proof.}
Note that the composition
$
H^i(X,\calO_X)\to\Nyg{}\Hcris{i+1}(X)\xto{}{F/p-1}\Hcris{i+1}(X)
$
is equal to the composition $H^i(X,\calO_X)\to\Nyg{}\Hcris{i+1}(X)\xto{}{F/p}\Hcris{i+1}(X),$ as can be seen from the long exact sequence \eqref{sel2}. Thus we must compare the map $\overline{F}$ with the restriction of $F/p$ to the image of $H^i(X,\calO_X)\to\Nyg{}\Hcris{i+1}(X)$.

Recall from Construction \ref{costruzionsfsup} that we first defined $F/p$ for eqrsp algebras, and then deduced its value for smooth schemes by descent. Thus, to compare $\overline{F}$ and $F/p$, we must make sense of $\overline{F}$ for quasisyntomic schemes. 

If $C$ is eqrsp, we see from tensoring the exact sequence
\[
0\to\Nyg{}\Acr(C)\to\Acr(C)\to C\to0
\]
with $\Z/p$ that $\Nyg{}\Acr(C)/p$ contains canonically a copy of $C$. Let $i:C\to\Nyg{}\Acr(C)/p$ denote this inclusion. In the notation of Proposition \ref{acrisspqr}, $i$ is $B$-linear and maps $x^{\alpha}$ to $px^{\alpha}$. This map is canonical and functorial for maps between eqrsp algebras. Thus we may proceed as in Construction \ref{costruzionsfsup} and define a canonical map $i:\Rgam(X,\calO_X)\to\Nyg{}\Rcris(X)\otimes^L\Z/p$ for any smooth or eqrsp $X$. Now consider the diagram
\[\begin{tikzcd}
	{\Rgam(X,\calO_X)} \\
	{\Nyg{}\Rcris(X)\otimes^L\Z/p} & {\Nyg{}\Rcris(X)[1]} & {\Nyg{}\Rcris(X)[1]} \\
	{\Rcris(X)\otimes^L\Z/p} & {\Rcris(X)[1]} & {\Rcris(X)[1]}
	\arrow["i"', from=1-1, to=2-1]
	\arrow[from=1-1, to=2-2]
	\arrow[from=2-1, to=2-2]
	\arrow["{F/p}"', from=2-1, to=3-1]
	\arrow["{\cdot p}", from=2-2, to=2-3]
	\arrow["{F/p}"', from=2-2, to=3-2]
	\arrow[from=2-3, to=3-3]
	\arrow[from=3-1, to=3-2]
	\arrow["{\cdot p}", from=3-2, to=3-3]
\end{tikzcd}\]
where the diagonal map comes from \eqref{triangoloagain}. By definition of $i$ it commutes when $X$ is eqrsp, so it also commutes for smooth $X$. If we let $\tilde{F}$ denote the composition of the two leftmost vertical maps, we are left to prove that for smooth $X$, the maps $\tilde{F}$ and $\overline{F}$ defined above coincide on cohomology. First we claim that the diagram
\[\begin{tikzcd}
	{\Rcris(X)\otimes^L\Z/p} && {\Rcris(X)\otimes^L\Z/p} \\
	& {\Rgam(X,\calO_X)}
	\arrow["F", from=1-1, to=1-3]
	\arrow[from=1-1, to=2-2]
	\arrow["{\tilde{F}}"', from=2-2, to=1-3]
\end{tikzcd}\]
commutes, where the unnamed map is the usual augmentation map. By descent it suffices to check this for an eqrsp algebra $C$. Using the description of $i$ given above, we see that $\tilde{F}:C\to\Acr(C)/p$ is Frobenius-linear and sends $x^{\alpha}$ to $x^{p\alpha}$, in the notation of Proposition \ref{acrisspqr}. In particular, it is a ring homomorphism. The projection $\Acr(C)/p\to C$, on the other hand, is $B$-linear, and maps $\frac{x^{\alpha}}{(\alpha!)_p}$ to $0$ if $\alpha\ge1$ and to $x^{\alpha}$ otherwise. Composing we get the Frobenius map described in Proposition \ref{acrisspqr}, so the claim is proved.

Now let $Y$ be smooth. We have $\Rcris(Y)\otimes^L\Z/p\simeq\RdR(Y/k)$, and by Zariski descent we have $\overline{F}=\tilde{F}$ on cohomology for $Y$ provided we know it for any affine smooth scheme. So suppose that $Y$ is affine. Then, looking at $H^0$, which is the only relevant group, we get a commutative diagram
\[\begin{tikzcd}
	{\HdR{0}(Y/k)=R^p} && {R^p=\HdR{0}(Y/k)} \\
	& R
	\arrow["{x\mapsto x^p}", from=1-1, to=1-3]
	\arrow[hook, from=1-1, to=2-2]
	\arrow["{H^0(\tilde{F})}"', from=2-2, to=1-3]
\end{tikzcd}\]
Since $\tilde{F}:C\to\Acr(C)/p$ is a ring homomorphism when $C$ is eqrsp, by descent we have that $H^0(\tilde{F})$ is also a ring homomorphism. The diagram tells us that $H^0(\tilde{F})(x)^p=H^0(\tilde{F})(x^p)=x^{p^2}$ for all $x\in R$, and since $R$ is reduced we have $H^0(\tilde{F})(x)=x^p$ for all $x\in R$. Thus $H^0(\tilde{F})=H^0(\overline{F})$ and we are done.
$\qed$

\begin{prop}\label{finitedegreetwo}

The group $\Hflat{2}(X,\Zp(1))$ is a $\Zp$-module of finite type.

\end{prop}

\noindent{\it Proof.} We first prove that 
\[
F/p-1:\Nyg{}\Hcris{1}(X)\to\Hcris{1}(X)
\]
is surjective. Recall that $\Hcris{1}(X)$ is a Dieudonné module, and that by Lemma \ref{nygverschiebung} if $x\in\Hcris{1}(X)$ then $Vx\in\Nyg{}\Hcris{1}(X)$ and $(F/p)(Vx)=x$. Define $M_{<1}\subseteq\Hcris{1}(X)$ to be the set of elements on which $V$ is topologically nilpotent. Then if $a\in M_{<1}$ we have
\[
a=(F/p-1)(Va+V^2a+\dots),
\]
and the argument lies in $\Nyg{}\Hcris{1}(X)$ by Lemma \ref{nygverschiebung}, so $M_{<1}$ is in the image of $F/p-1$. Denote by $M_1$ the quotient $\Hcris{1}(X)/{M_{<1}}$. Since the sequence $(V^m(p^na))_m$ converges to $0$ if and only if the sequence $(V^m(a))_m$ converges to $0$, this is a torsion-free $\W(k)$-module, and it comes with the induced maps $\overline{F}$ and $\overline{V}$. If we show that
\begin{equation}\label{fsup}
\overline{F}/p-1:\overline{V}(M_1)\to M_1
\end{equation}
is surjective, we are done. But by design $M_1$ has no non-trivial $\overline{V}$-topologically nilpotent elements. So by the Dieudonné-Manin classification it is pure of slope $1$, with a basis $a_i$ such that $\overline{F}(a_i)=pa_i$. If $\lambda\in\W(k)$ we have $(\overline{F}/p-1)(\overline{V}(\lambda a_i))=(\lambda-\sigma^{-1}(\lambda))a_i$. The map $\W(k)\to \W(k)$ which maps $\lambda$ to $\lambda-\sigma^{-1}(\lambda)$ is surjective, so the map \eqref{fsup} is indeed surjective.

It follows that 
\[
0\to\Hflat{2}(X,\Zp(1))\to\Nyg{}\Hcris{2}(X)\xto{}{F/p-1}\Hcris{2}(X)
\]
is exact. We have an exact sequence
\[
0\to L\to\Nyg{}\Hcris{2}(X)\to\Hcris{2}(X)
\]
where $L=\coker\left(q_1:\Hcris{1}(X)\to H^1(X,\calO_X)\right)$. Suppose that $F/p-1$ is injective on $L$, or equivalently that $\ker(\Nyg{\Hcris{2}(X)\to\Hcris{2}(X))\cap\ker(F/p-1)=0}$. Then the composition
\begin{equation*}
    \Hflat{2}(X,\Zp(1))\to\Nyg\Hcris{2}(X)\to \Hcris{2}(X)
\end{equation*}
is injective. Its image is contained in $\Hcris{2}(X)^{F=p}$, so we have
\[
\Hflat{2}(X,\Zp(1))\subseteq\Hcris{2}(X)^{F=p},
\]
and since the right-hand $\Zp$-module is of finite type by the Dieudonné-Manin classification we are done. So let us prove that $F/p-1$ is injective on $L$.

Let $a\in\ker(F/p-1|_L)$ and lift it to $a'\in H^1(X,\calO_X)$. Lemma \ref{fsupcoerente} implies that $\overline{F}(a)$ lies in $\ker\left(\HdR{1}(X/k)\to\Hcris{2}(X)[p]\right)$,
and by Proposition \ref{confrontoderham} we have
\[
\ker\left(\HdR{1}(X/k)\to\Hcris{2}(X)[p]\right)=\im\left(q:\Hcris{1}(X)\to\HdR{1}(X/k)\right),
\]
so we may write $\overline{F}(a')=q(b)$ for some $b\in\Hcris{1}(X)$. Going back to the commutative diagram \ref{diagrammacrisderham}, we see that  $V(b)$ maps to zero in $\HdR{1}(X/k)$, so $V(b)=pc$ for some $c\in\Hcris{1}(X)$, hence $pb=pF(c)$ and finally $b=F(c)$ since $\Hcris{1}(X)$ is torsion-free. Thus $\overline{F}(a')=q(F(c))=\overline{F}(q_1(c))$ and by injectivity of $\overline{F}$ we have $a'=q_1(c)$, which implies $a=0$. This is what we wanted to prove.
$\qed$

With a little more effort the proof of Proposition \ref{finitedegreetwo} yields the following result of Illusie, see \cite[II.5.14]{illusiedrw}.

\begin{prop}\label{}

If $\Hcris{2}(X)$ is torsion-free and every global differential $1$-forms is closed, there is a short exact sequence
\[
0\to\Hflat{2}(X,\Zp(1))\to\Hcris{2}(X)\xto{}{F-p}\Hcris{2}(X).
\]
where the first map is the natural comparison map.
\end{prop}

\noindent{\it Proof.} Keep the notation from the proof of Proposition \ref{finitedegreetwo}. First we claim that $\Hcris{1}(X)\to H^1(X,\calO_X)$ is surjective. Since $\Hcris{2}(X)[p]=0$ the composition \eqref{comp2} for $i=1$ is zero, so by Lemma \ref{fsupcoerente} we have that $F/p-1$ is zero on $L=\coker\left(\Hcris{1}(X)\to H^1(X,\calO_X)\right)$. But we have proved that $F/p-1$ is injective on $L$, so $L=0$.

Thus $\Nyg{}\Hcris{2}(X)\subseteq\Hcris{2}(X)$ and both groups are torsion-free. We have a commutative diagram
\[\begin{tikzcd}
	0 & {\Hflat{2}(X,\Zp(1))} & {\Nyg{}\Hcris{2}(X)} & {\Hcris{2}(X)} \\
	0 & K & {\Hcris{2}(X)} & {\Hcris{2}(X)} \\
	&& {H^2(X,\calO_X)} & {\HdR{2}(X/k)}
	\arrow[from=1-1, to=1-2]
	\arrow[from=1-2, to=1-3]
	\arrow["f"', from=1-2, to=2-2]
	\arrow["{F/p-1}", from=1-3, to=1-4]
	\arrow[hook, from=1-3, to=2-3]
	\arrow["{\cdot p}"', from=1-4, to=2-4]
	\arrow[from=2-1, to=2-2]
	\arrow[from=2-2, to=2-3]
	\arrow["{F-p}", from=2-3, to=2-4]
	\arrow[from=2-3, to=3-3]
	\arrow[from=2-4, to=3-4]
	\arrow["{\overline{F}}", from=3-3, to=3-4]
\end{tikzcd}\]
with exact rows, with the two right-hand columns exact, and where all vertical maps are injective. We want to show that $f$ is surjective and an easy diagram chase shows that it is enough to show that $\overline{F}$ is injective. This is true by Lemma \ref{comportamentofbar}.

$\qed$

\subsection{Straight varieties and abelian varieties}\label{straight} Now we want to single out a class of varieties for which the Nygaard filtration and fppf cohomology are determined uniquely by crystalline cohomology. In general these groups depend on crystalline cohomology and its relation with de Rham cohomology, which can be quite subtle. We call these varieties straight, but they have appeared in previous literature as Mazur--Ogus varieties \cite{ekedahl}. Toward the end of this section we show how to explicitly compute the fppf cohomology of some very simple straight varieties, namely ordinary abelian varieties and products of supersingular elliptic curves. Here is the definition.

\begin{deffa}\label{}

A variety $X$ is called \textit{straight} if the following conditions are satisfied:
\begin{enumerate}
    \item For all $i$ the group $\Hcris{i}(X)$ is torsion-free.
    \item The Hodge-de Rham spectral sequence for $X$ degenerates at page $E_1$.
\end{enumerate}

\end{deffa}

By the discussion preceding Lemma \ref{bordoconiugata}, the conjugate spectral sequence of a straight variety also degenerates. Example of straight varieties are abelian varieties - see \cite[Proposition 5.1]{oda} for the degeneration of the Hodge-de Rham spectral sequence, and see below for the claim about crystalline cohomology - and K3 surfaces \cite[II.7.2]{illusiedrw}. 

Let $X$ be a straight variety. Condition (1) and Proposition \ref{confrontoderham} imply that the maps $\Hcris{i}(X)\to\HdR{i}(X/k)$ are surjective, while condition (2) ensures that $\HdR{i}(X/k)\to H^{i}(X,\calO_X)$ is surjective. From the long exact sequence \eqref{sel2} it follows that
\begin{equation}\label{nygker}
\Nyg{}\Hcris{i}(X)=\ker\left(\Hcris{i}(X)\to H^i(X,\calO_X)\right)
\end{equation}
for all $i$, so $\Nyg\Hcris{i}(X)$ is torsion-free as well. In Section \ref{pruffa} we defined maps $F/p:\Nyg{}\Hcris{i}(X)\to\Hcris{i}(X)$. Their existence implies that $\Nyg{}\Hcris{i}(X)$ is contained in
\[
F^{-1}(p\Hcris{i}(X))=\{a\in\Hcris{i}(X)\text{ s.t. }p|F(a)\}
\]
for all $i$, but in fact the reverse inclusion also holds true.

\begin{prop}\label{nygstraight}
If $X$ is a straight variety, the equality
\[
\Nyg{}\Hcris{i}(X)=F^{-1}(p\Hcris{i}(X))=\{a\in\Hcris{i}(X)\text{ s.t. }p|F(a)\}
\]
holds for all $i$.
\end{prop}

\noindent{\it Proof.}
Taking cohomology in the commutative diagram
\[
\begin{tikzcd}
	{\Nyg{}\Rcris(X)} & {\Rcris(X)} & {\Rgam(X,\calO_X)} \\
	{\Rcris(X)} & {\Rcris(X)} & {\RdR(X/k)}
	\arrow[from=1-1, to=1-2]
	\arrow["{F/p}", from=1-1, to=2-1]
	\arrow[from=1-2, to=1-3]
	\arrow["F", from=1-2, to=2-2]
	\arrow["{\overline{F}}", from=1-3, to=2-3]
	\arrow["{\cdot p}", from=2-1, to=2-2]
	\arrow[from=2-2, to=2-3]
\end{tikzcd}
\]
whose rows are distinguished, we get the following commutative diagram with exact rows:
\[
\begin{tikzcd}
	0 & {\Nyg{}\Hcris{i}(X)} & {\Hcris{i}(X)} & {H^{i}(X,\calO_X)} & 0 \\
	0 & {\Hcris{i}(X)} & {\Hcris{i}(X)} & {\HdR{i}(X/k)} & 0
	\arrow[from=1-1, to=1-2]
	\arrow[from=1-2, to=1-3]
	\arrow["{F/p}", from=1-2, to=2-2]
	\arrow[from=1-3, to=1-4]
	\arrow["F", from=1-3, to=2-3]
	\arrow[from=1-4, to=1-5]
	\arrow["{\overline{F}}", from=1-4, to=2-4]
	\arrow[from=2-1, to=2-2]
	\arrow["{\cdot p}", from=2-2, to=2-3]
	\arrow[from=2-3, to=2-4]
	\arrow[from=2-4, to=2-5]
\end{tikzcd}
\]
Here we have used \eqref{nygker} for the top row and Proposition \ref{confrontoderham} for the bottom row We thus have a map
\[
\Nyg{}\Hcris{i}(X)\xto{}{g=(F/p,1)}\{(a,b)\in\Hcris{i}(X)^{\oplus2}\text{ s.t. }F(b)=pa\},
\]
and we must prove it is an isomorphism. A pair $(a,b)$ in the right-hand group is uniquely determined by $b$. Thus $g$ is injective. Viceversa, if $(a,b)$ is such a pair, by injectivity of $\overline{F}$, which follows from the degeneration of the conjugate spectral sequence by Lemma \ref{bordoconiugata}, we see that $b$ must come from $\Nyg{}\Hcris{i}(X)$. Therefore $g$ is also surjective.
$\qed$

In other words the Nygaard filtration is uniquely determined by the $F$-crystal $\Hcris{i}(X)$ for a straight variety $X$. Thus fppf cohomology with $\Zp(1)$ coefficients is uniquely determined by the $F$-crystal $\Hcris{i}(X)$. To be precise, we have the following statement.

\begin{lem}\label{fppfstraight}

Let $X$ be a straight variety. For all $i\ge0$ there is a short exact sequence
\[
0\to\Hflat{i}(X,\Zp(1))/\text{tors}\to\Nyg{}\Hcris{i}(X)\xto{}{F/p-1}\Hcris{i}(X)\to\Hflat{i+1}(X,\Zp(1))_{\text{tors}}\to0.
\]
In particular the fppf cohomology of $X$ only depends on the crystalline cohomology of $X$ with its Frobenius endomorphism.

\end{lem}

\noindent{\it Proof.}
This follows by combining \eqref{sel1} with the fact that $\Nyg{}\Hcris{1}(X)$ is torsion-free, and with Lemma \ref{lemma2} which says that $\coker\left(F/p-1\right)$ is torsion.
$\qed$

We now turn our attention to abelian varieties. If $A$ is an abelian variety we have $\Hcris{i}(A)\simeq\wedge^i\Hcris{1}(A)$ (\cite[II.7.1]{illusiedrw}), and $\Hcris{1}(A)$ is the contravariant Dieudonné module of $A[p^\infty]$ by Proposition \ref{crisabeliana}. Therefore the cohomology groups $\Hflat{i}(A,\Zp(1))$ only depend on the Dieudonné module of $A$.

\begin{cor}\label{moltiplicazionen}

Let $A$ be an abelian variety over $k$. The multiplication-by-$n$ map $[n]:A\to A$ acts as $n^{i}$ on $\Hflat{i+1}(A,\Zp(1))_{\text{tors}}$ and on $\Hflat{i}(A,\Zp(1))/\text{tors}$.

\end{cor}

\noindent{\it Proof.}
If we know that $[n]$ acts as $n^i$ on $\Hcris{i}(A)$ and on $\Nyg{}\Hcris{i}(A)$ the conclusion follows from Lemma \ref{fppfstraight}. Since $\Hcris{i}(A)\simeq\wedge^i\Hcris{1}(A)$ it is enough to check this for $i=1$, where it is a consequence of the fact that Dieudonné theory is an equivalence of categories. By Proposition \ref{nygstraight} $\Nyg{}\Hcris{i}(X)$ is a subgroup of $\Hcris{i}(A)$, thus $[n]$ also acts as $n^i$ on $\Nyg{}\Hcris{i}(X)$.
$\qed$

In a future paper with A.Skorobogatov and Y.Yang we will use these techniques to study the $p$-torsion of the Brauer group of an abelian variety in more detail. We end this section with some computations of fppf cohomology in two simple cases.

\vspace{.5em}
\noindent\textbf{First computation.} Let $A$ be an ordinary abelian variety of dimension $g$ over $k$, i.e. $A[p](k)\simeq(\Z/p)^g$. It is well-known that $\Hcris{1}(A)$ has a basis $(x_1,\dots,x_g,y_1,\dots,y_g)$ such that $F(x_i)=x_i$ and $F(y_i)=py_i$. We compute the fppf cohomology of $A$.

\begin{prop}\label{ordinaria}

The following hold:
\begin{itemize}
    \item[(1)] The $\W(k)$-module $\Hcris{i}(A/\W(k))$ is free with basis the $x_I\wedge y_J$, where $I,J\subseteq\{1,\dots,g\}$ and $|I|+|J|=i$.
    \item[(2)] The $\W(k)$-module $\Nyg\Hcris{i}(A/\W(k))$ is free with basis the $p^{\varepsilon}x_I\wedge y_J$, where $I,J\subseteq\{1,\dots,g\}$, $|I|+|J|=i$, and $\varepsilon=1$ if $J=\emptyset$, $\varepsilon=0$ otherwise.
    \item[(3)] The $\Zp$-module $\Hflat{i}(A,\Zp(1))$ is free of rank $g\cdot\binom{g}{i-1}$.
\end{itemize}

\end{prop}

\noindent{\it Proof.} Statement (1) follows from the identification $\Hcris{i}(A)\simeq\wedge^i\Hcris{1}(A)$. Then (2) follows from (1) and Proposition \ref{nygstraight}. 

For (3), we compute $F/p-1$ on the system of generators of $\Nyg{}\Hcris{i}(X)$ we have found in (2). If $\lambda\in\W$,
\[
(F/p-1)(\lambda p^{\varepsilon}x_I\wedge y_J)=(p^{\varepsilon+|J|-1}\sigma(\lambda)-\lambda)x_I\wedge y_J.
\]
As $\lambda$ varies in $\W$, $(p^{\varepsilon+|J|-1}\sigma(\lambda)-\lambda)$ covers the whole of $\W$, so $F/p-1$ is surjective. 

On the other hand $\ker\left(F/p-1\right)$ is a free $\Zp$-module with basis the $x_I\wedge y_J$ where $|J|=1$. From Lemma \ref{fppfstraight} we get (3).
$\qed$

\vspace{.5em}
\noindent\textbf{Second computation.} Let $E$ be a supersingular elliptic curve over $k$, i.e. $\Hcris{1}(E)$ is pure of slope $1/2$, and set $A=E\times_kE$. We will show that $\Hflat{3}(A,\Zp(1))$ is an infinite $p$-torsion group which can be identified with $k$. Referring back to Remark \ref{unipotente} this means that $U_3\simeq\Ga$. More general computations for abelian surfaces can be found in \cite[II.7.1]{illusiedrw}.

A simple argument shows that $\Hcris{1}(E/\W(k))$ has a $\W(k)$-basis $(x,y)$ such that $F(x)=y$ and $F(y)=px$. Indeed we know that $\Hcris{1}(E/\W(k))[1/p]$ has a $\K$-basis $(a,b)$ such that $F(a)=b$ and $F(b)=pa$. We may take $a\in\Hcris{1}(E/\W(k))$ not divisible by $p$, and then $b=F(a)$ also lies in $\Hcris{1}(E/\W(k))$. If $b$ is not divisible by $p$, set $(x,y)=(a,b)$. Otherwise set $(x,y)=(b/p,a)$. By construction neither $x$ nor $y$ is divisible by $p$ in $\Hcris{1}(E/\W(k))$, which is a free $\W(k)$-module of rank $2$, thus $(x,y)$ has the desired property.

By the crystalline Künneth formula \cite[V.4.2.3]{berthelotcris}, the $\W(k)$-module $\Hcris{2}(A/\W(k))$ has a basis $(a_1,a_2,b_1,b_2,c,d)$ such that
\begin{enumerate}
    \item $F(a_1)=p^2a_1,F(a_2)=p^2a_2$
    \item $F(b_1)=b_2,F(b_2)=p^2b_1$
    \item $F(c)=pc,F(d)=-pd$
\end{enumerate}
We may take for example $b_1=x\otimes x,b_2=y\otimes y$, then $c=(x\otimes y+y\otimes x)/2$ and $d=(x\otimes y-y\otimes x)/2$. Then by Proposition \ref{nygstraight} the $\W(k)$-module $\Nyg\Hcris{2}(A/\W(k))$ is generated by $a_1,a_2,pb_1,b_2,c,d$.

\begin{prop}\label{gruppoinifinito}

The group $\Hflat{3}(A,\Zp(1))$ is torsion, infinite and isomorphic to $k$.

\end{prop}

\noindent{\it Proof.} The rank of $\Hflat{3}(A,\Zp(1))$ is $0$ because $\Hcris{3}(A)$ is pure of slope $3/2$. Therefore $\Hflat{3}(A,\Zp(1))$ is torsion and we have
\[
\Hflat{3}(A,\Zp(1))_{\text{tors}}\simeq\coker\left(F/p-1:\Nyg\Hcris{2}(A/\W(k))\to\Hcris{2}(A/\W(k))\right)
\]
Let's determine the image of $f=F/p-1$: if $\lambda\in\W(k)$ we have:
\begin{gather*}
f(\lambda a_1)=a_1(p\sigma(\lambda)-\lambda)\\
f(\lambda a_2)=a_2(p\sigma(\lambda)-\lambda)\\
f(\lambda pb_1)=\sigma(\lambda)b_2-\lambda pb_1\\
f(\lambda b_2)=\sigma(\lambda)pb_1-\lambda b_2\\
f(\lambda c)=c(\sigma(\lambda)-\lambda)\\
f(\lambda d)=-d(\sigma(\lambda)-\lambda)
\end{gather*}
Thus the image of $f$ is contained in $\Nyg\Hcris{2}(A/\Zp)$ and it is easy to see that as $\lambda$ varies we get a set of generators of this group. Therefore $\coker\left(F/p-1\right)$ is isomorphic to $b_1\W(k)/p\W(k)\simeq k$.
$\qed$

Using the same steps it is of course possible to compute the other cohomology groups $\Hflat{i}(A,\Zp(1))$, which turn out to be torsion-free for $i\ne3$. The analogous computations for arbitrary products of supersingular elliptic curves are very similar. We should mention that in \cite{yuanbrauer} the author computes $\Hflat{3}(X,\Zp(1))_{\text{tors}}$ for these varieties and many other using different techniques.

\vspace{.5em}
\hrule
\vspace{1em}
\eu

\end{document}